\documentclass{elsarticle}
\usepackage{graphicx}
\usepackage{subfigure}
\usepackage{color}
\usepackage{newlfont}
\usepackage{multirow}
\usepackage{longtable} 
\usepackage{ifthen}
\usepackage{alltt}
\usepackage{enumerate}
 \usepackage[text={6.25in,8.5in},centering]{geometry} 

\newcommand{\ba}{\begin{array} }
\newcommand{\ea}{\end{array} }
\newcommand{\bae}{\begin{eqnarray}}
\newcommand{\eae}{\end{eqnarray}}
\newcommand{\bea}{\begin{eqnarray*}}
\newcommand{\eea}{\end{eqnarray*}}
\newcommand{\be}{\begin{equation}}
\newcommand{\ee}{\end{equation}}

\newcommand{\modifyb}[1]{\textcolor{black}{#1}}
\newcommand{\modifyr}[1]{\textcolor{black}{#1}}
\newcommand{\modifyc}[1]{\textcolor{black}{#1}}

\newcommand{\pr}{{\bf Proof}~~}

\usepackage{amssymb}
\usepackage{amsthm}
\usepackage{amsmath}

\usepackage{graphicx}
\usepackage{color}
\usepackage[colorlinks]{hyperref}
\usepackage{lscape}
\usepackage{graphicx}
\usepackage{epstopdf}
\newtheorem{theorem}{\hskip\parindent\bf Theorem}[section]
\newtheorem{lemma}{\hskip\parindent\bf Lemma}[section]
\newtheorem{proposition}{\bf Proposition}[section]

\newtheorem{corollary}{\hskip\parindent\bf Corollary}[section]

\begin{document}
 \markboth{Allee effects in two species interaction models}{}
 \title{Dynamics of \modifyb{a generalized Beverton-Holt competition model} subject to Allee effects}
\author{Yun Kang\footnote{Sciences and Mathematics Faculty, Arizona State University, Mesa, AZ 85212, USA. E-mail: yun.kang@asu.edu}}
\begin{abstract}
In this article, we propose and study \modifyb{a generalized Beverton-Holt competition model} subject to Allee effects to obtain insights on how the interplay of Allee effects and contest competition affects the persistence and the extinction of two competing species. By using \modifyr{the theory of monotone dynamics} and the properties of critical curves for non-invertible maps, our analysis shows that our model has relatively simple dynamics, i.e., almost every trajectory converges to a \modifyr{locally asymptotically stable} equilibrium if \modifyr{the intensity of intra-specific competition intensity exceeds that of inter-specific competition. This equilibrium dynamics is also possible when the intensity of intra-specific competition intensity is less than that of inter-specific competition but under conditions that  the maximum intrinsic growth rate of one species is not too large. The coexistence of two competing species occurs only if the system has four interior equilibria.} We provide an \modifyr{approximation to} the basins of the boundary attractors (i.e., the extinction of one or both species) where our results suggests that  \modifyr{contest species are more prone to extinction than scramble ones are at low densities}. In addition, \modifyr{in comparison to the} dynamics of two species scramble competition models subject to Allee effects, our study suggests that (i) Both contest and scramble competition models can have only three boundary attractors without the coexistence equilibria, or four attractors among which only one is the \modifyr{persistent attractor}, \modifyr{whereas} scramble competition \modifyr{models} may have the extinction of both species as its only attractor under certain conditions, i.e., \emph{the essential extinction} of two species due to strong Allee effects; (ii) \modifyr{Scramble competition models like Ricker type models can} have much more complicated dynamical structure of interior attractors than contest ones like Beverton-Holt type models have; 
and (iii) Scramble competition models \modifyb{like Ricker type competition models} may be more likely to promote the coexistence of two species at low and high densities \modifyc{under certain conditions: At low densities, weak Allee effects decrease the fitness of resident species so that the other species is able to invade at its low densities; While at high densities, scramble competition can bring the current high population density to a lower population density but is above the Allee threshold in the next season, which may rescue a species that has essential extinction caused by strong Allee effects.} Our results may have potential to be useful for conservation biology: For example, \modifyc{if one endangered species is facing essential extinction due to strong Allee effects, then we may rescue this species by bringing another competing species subject to scramble competition and Allee effects under certain conditions.}


 \end{abstract}

\bigskip
\maketitle
\section{Introduction}

An Allee effect is a biological phenomenon characterized by a positive correlation between a population density and its per capita growth rate at small population densities (Allee \emph{et al.} 1949). A distinction is made between a \emph{strong Allee effect} and a \emph{weak Allee effect}:  a \emph{strong Allee effect} refers as to a population that exhibits a ``critical size or density", below which population declines to extinction, and above which it may increase; while a \emph{weak Allee effect} refers as to a population that lacks a ``critical density", but where, at lower densities, the population growth rate rises with increasing densities (Stephens \emph{et al.} 1999; Lidicker 2010).  Allee effects have been detected in natural populations for a wide array of taxa and are believed to be strong regulators of extinction, influencing colonization success of invasive species, disease dynamics, and the long-run population demise in conservation. Mathematical models can help us understand these combined effects on persistence, expansion or extinction of biological species. There \modifyr{is a considerable amount of} populations dynamics in the presence of Allee effects (e.g., Dennis 1989 \& 2002; Selgrade and Namkoong 1992; McCarthy 1997; Shigesada and Kawasaki 1997; Greene and Stamps 2001; Keitt \emph{et al} 2001; Fagan \emph{et al} 2002; Wang \emph{et al} 2002; Liebhold and Bascompte 2003; Schreiber 2003; Drake 2004; Zhou \emph{et al} 2004; Petrovskii \emph{et al} 2005; Taylor and Hastings 2005; Jang 2006; Aguirre \emph{et al} 2009; Egami 2009\&2010; Thieme \emph{et al.} 2009; Elaydi and Sacker 2010; Wang \emph{et al} 2010; Liu \emph{et al} 2011; Kang and Castillo-Chavez 2012; Kang and Castillo-Chavez 2014a\&b; Cai \emph{et al} 2014; Peng and Kang 2015; Livadiotis \emph{et al }2015) as well as various models in patchy environments (e.g., Amarasekare 1998a \&1998b; Gyllenberg \emph{et al} 1999; Ackleh \emph{et al} 2007; Kang and Lanchier 2011; Kang and Armbruster 2011; Kang  and Castillo-Chavez 2014c). In another direction, evolutionary models with Allee effects have gained some attentions (Cushing and Hudson 2012; Kang and Udiani 2014; Cushing 2015). For example,  Cushing and Hudson (2012)  investigated the global dynamics of an evolutionary model for a population subject to strong Allee effects. One of their results suggests that evolution is beneficial in the sense that reduces the possibility of extinction due to an Allee effect. \\

Resource competition can be defined in the contest-scramble spectrum \modifyr{competition} according to resource partitioning among competitors (Nicholson 1954). Contest and scramble-type \modifyr{competition} are two extreme forms of competition, which are characterized by resource monopolization and resource sharing, respectively. Contest competition results in a constant number of survivors with enough resource gain against an initial density of competitors, while scramble competition (also called over-compensation) results in an increased number of survivors with decreased resource intake and an increased number of competitors (Calow \emph{et al.} 1998).\modifyc{ A well-known example for contest competition is the Beverton-Holt model while the well-known example for scramble competition is the Ricker model. Even though Allee effects are expected commonly in nature, they are often ignored in studies of species coexistence mechanisms.  Many species can experience \modifyr{different forms of competition at high densities and Allee effects at low densities with the consequence} that their populations do not grow optimally at low densities and individuals compete with one another at high densities (Begon \emph{et al.} 1996; Etiemme \emph{et al.} 2002; Kang 2013). The potential effects of combining positive density dependences from Allee effects with negative density dependences coming from different types of competition are expected to generate rich dynamics due to the fact that outcomes of varied patterns of competition have different consequences for population dynamics and the evolution of individual traits (Godfray 1987; Lomnicki 1988; Ives 1989).  For instance, with a species at a low density, in the invader state, it faces both intra-specific and inter-specific competition, but obtains a lesser effect of reduced intra-specific competition due to Allee effects compared with the standard investigations of models of competing species. As comparisons between intra-specific and inter-specific density dependences are critical to competitive coexistence, ignoring the potential that Allee effects greatly modify intra-specific competition is a critical oversight.  For instance, Chesson and Ellner (1989) has pointed out the potential for major effects such as stochastic models may have interior invariant regions if low-density growth rates are depressed by Allee effects. Hopf \emph{et al} (1993) also point out potential major implications for community structure of Allee effects uncovered by introducing them into a competition model. }\\

\modifyc{Recently, Kang and Yakubu (2011) and Kang (2013) have explored the population dynamics of two competing species when both species experience \emph{scramble competition} (i.e., Ricker-type models) and \emph{Allee effects} induced by predation \modifyr{saturations}. The positive density dependence from weak Allee effects can decrease the fitness of resident species so that the other species is able to invade at its low densities. Thus \emph{weak Allee effects} can promote the permanence of two competing species at their low densities (Kang and Yakubu 2011). In the presence of strong Allee effects, scramble inter-specific competition can bring the current hight population density to a lower population density but being above the Allee threshold. Thus \emph{strong Allee effects} can save species from potential \emph{essential extinction} (Schreiber 2003) by bringing in scramble inter-specific competition (Kang 2013), therefore promote the existence of two competing species at their high densities under certain conditions. Understanding how different \modifyr{forms of competition in the presence of Allee effects affect population dynamics} when species suffer from Allee effects, can advance our understanding of the extinction and establishment of species in ecological communities, with implications for conservation programs (Zhou \emph{et al} 2004; Courchamp \emph{et al.} 2009; Kang and Yakubu 2011; Kang 2013; Kang \emph{et al} 2014 a\&b). } \\

A recent work by Livadiotis and Elaydi (2012) considered a general framework of population models subject to strong Allee effects. \modifyr{We} adopt and modify their definitions on \modifyr{mathematical formulation of} strong Allee effects \modifyr{in single species models.} \modifyr{Our main purpose} is to investigate the \modifyr{dynamical outcomes in two species contest competition models} when both species have \emph{Allee effects}. \modifyb{The contest competition models in this article refer as to a generalized Beverton-Holt competition model while scramble competition models refer as to species that have one-hump growth functions (e.g., Ricker's map) in its single state.} More precisely, we would like to explore possible answers to the following ecological questions:
\begin{enumerate}
\item What are the dynamical outcomes of a two \modifyr{species competition model when both species} experience both strong Allee effects and contest \modifyr{competition}?
\item \modifyr{How do the dynamics of models with a strong Allee effect and contest competition compare to those of models with a strong Allen effect and scramble competition?}
\item What are the generic dynamical features of population models subject to strong Allee effects and \modifyr{different forms of competition at low and high densities?}\\

\end{enumerate}

\modifyc{The \modifyr{paper} is organized as follows: 
In Section 2, we \modifyr{define how we will mathematically formulate} Allee effects, \modifyr{contest and scramble competition} for single species population models and derive \modifyr{a generalized Beverton-Holt competition model} that each species is subject to strong Allee effects. In addition, we provide the important results of single species models with strong Allee effects and contest competition that can be applied to two species competition models. 
In Section 3, we apply \modifyr{the theory of monotone dynamics} to obtain sufficient conditions when our two species competition model has \modifyr{equilibrium} dynamics. This result can establish some generic dynamical features of contest competition models subject to Allee effects. We also give approximated basins of attractions of boundary attractors (i.e., initial conditions that lead to the extinction of one or both species). Moreover, we explore conditions that lead to the coexistence of two competing species. These results can help us understand the coexistence and the extinction conditions for two competing species that are subject to contest competition and Allee effects.   In Section 4, we apply the theoretical results derived in Section 3 to a specific two competing species model with strong Allee effects and contest \modifyr{competition}. Numerical simulations are performed to obtain more insights \modifyr{into} the generic features of two species models with strong Allee effects and contest \modifyr{competition}. The results obtained from Section 3 \&4 can partially answer three questions proposed in the introduction.
In Section 5, we discuss the dynamic outcomes of two competing species population models with different types of competition (i.e., contest and scramble competition) where each species has strong Allee effects. In addition, 
we \modifyr{summarize} our results and \modifyr{mention some} potential future studies. In Appendix, we provide \modifyr{detailed proofs of} our theoretical results.}
\section{Model Derivations}
In this article, we are interested in the population dynamics of a \modifyr{generalized Beverton-Holt two species competition model} subject to Allee effects. More specifically, we would like to explore how the interactions of contest competition and Allee effects affect the establishment and extinction of species. \modifyb{By understanding population dynamics of such models, we are able to obtain insights on the impact of different forms of competition combined with Allee effects on the persistence and extinction of species by comparing the dynamics to the case when two competing species both suffer from scramble competition and Allee effects}. 

\subsection{Single species population models}
\begin{figure}[ht]
\centering
\subfigure[]{
   \includegraphics[scale =0.395] {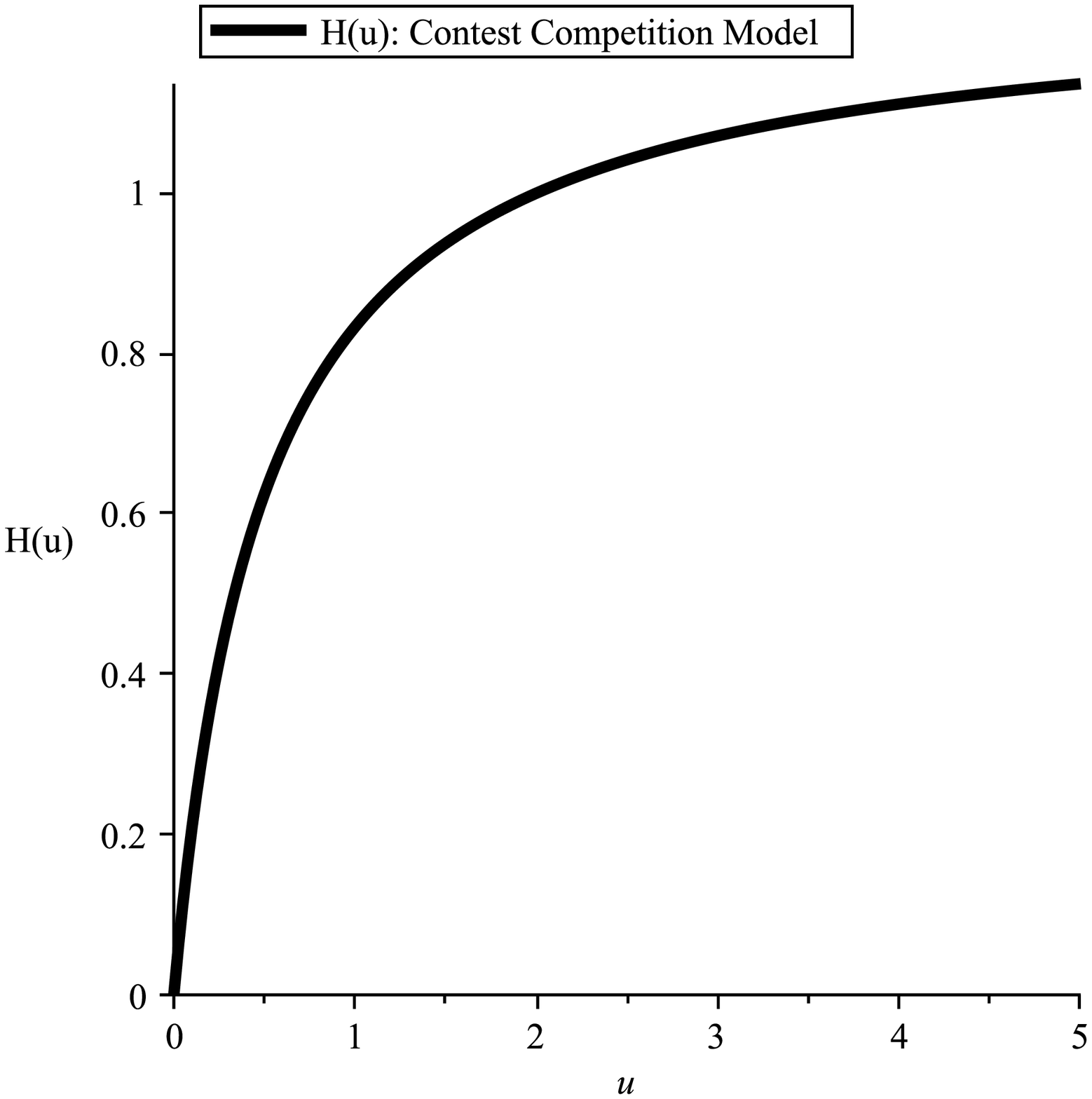} \label{fig1:contest}} 
   \subfigure[]{
   \includegraphics[scale =.395] {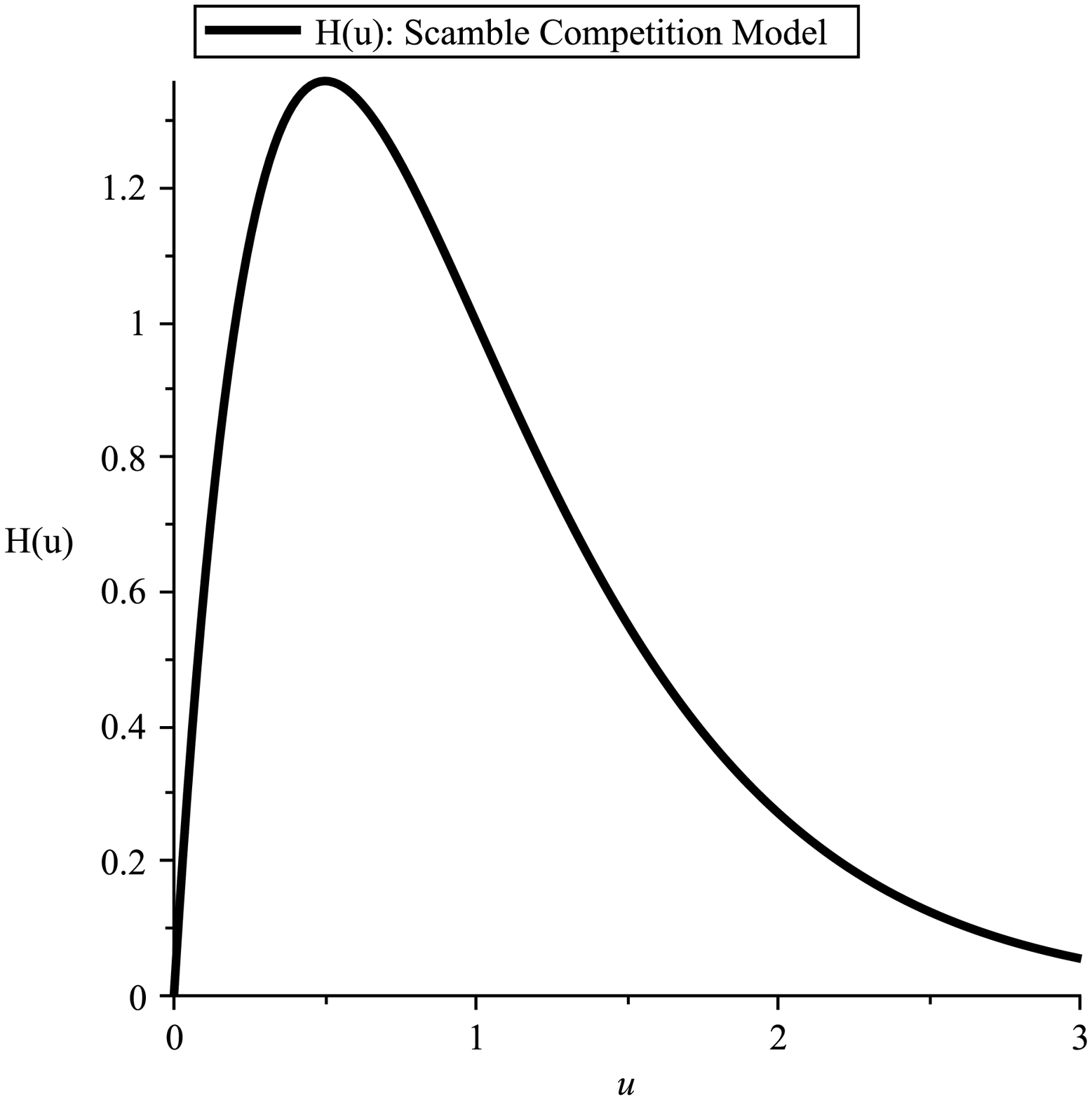}\label{fig1:scramble}}
  \caption{Examples of single species models with contest competition (a) versus scramble competition (b).} \label{fig1:com}
\end{figure}

First, we would like to introduce the terms of contest competition, scramble competition and Allee effects for one species population models rigorously. Let $u_t$ be a species population density in season $t$, then its population density in season $t+1$ can be represented as the following equation:
\bae\label{s}
u_{t+1}&=&H(u_t)=u_t h(u_t)
\eae where $H(u)\geq 0, u\in \mathbb R_+$ is the growth function of species $u$ and $h(u)$ is the per capital growth rate that is subject to the condition $\lim_{u\rightarrow}h(u)=0$ \modifyb{which guarantees the boundedness of population} (Kang and Chesson 2010; Kang 2013). We say that species $u$ has \emph{contest intra-specific competition} if $\lim_{u\rightarrow\infty} H(u)=c>0$ which \modifyb{implies that a constant number of survivors with enough resource gain against an initial density of competitors}. We say that species $u$ has \emph{scramble intra-specific competition} if $\lim_{u\rightarrow\infty} H(u)=0$ which \modifyb{implies that an increased number of survivors with decreased resource intake and an increased number of competitors.} \emph{Contest} and \emph{scramble} intra-specific competition are two extreme forms of competition that can result in very different consequences for population dynamics. For example, an important ecological example of contest intra-specific competition is the Beverton-Holt model (see Figure \ref{fig1:contest}) whose $H$ takes the form of $H(u)=\frac{r u}{1+au}$ while a well-known ecological example of scramble intra-specific competition is the Ricker's model (see Figure \ref{fig1:scramble}) whose $H$ takes the form of $H(u)=ue^{r-au}$. The Beverton-Holt model has relatively simple population dynamics where population goes to 0 if $r<1$ and population goes to $r-1$ if $r>1$. While the Ricker's model has very complicate population dynamics (e.g., periodic orbits, chaos) depending on the values of $r$. See Figure \ref{fig1:com} for examples of single species models with contest intra-specific competition versus scramble intra-specific competition.\\

\begin{figure}[ht]
\centering
   \includegraphics[scale =0.65] {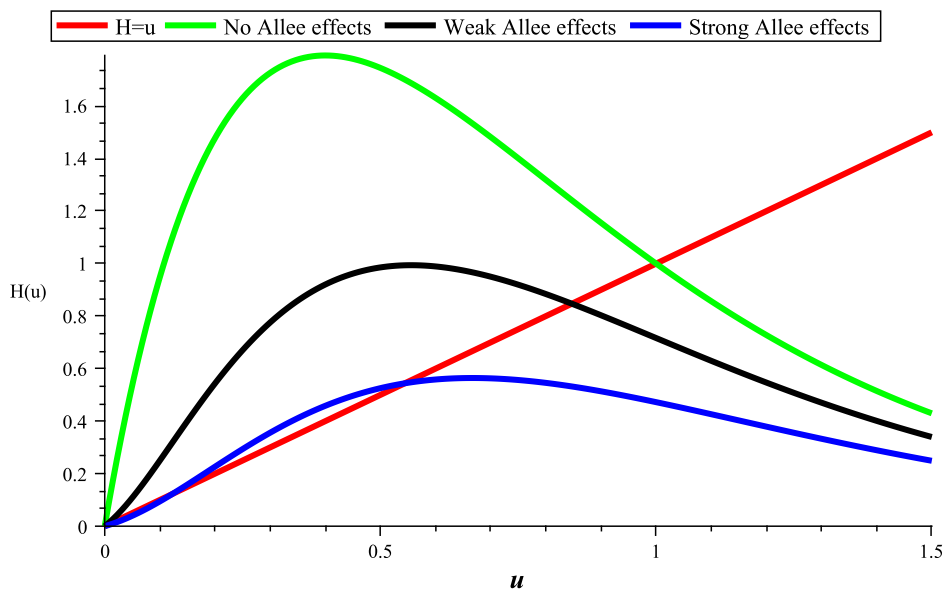}
   \caption{\modifyc{An illustration of strong Allee effects versus weak Allee effects of the growth function $H(u)=uh(u)$: The red curve is $H=u$ or $h=1$; the green curve is the case when \eqref{s} has no Allee effects; the blue curve is the case when \eqref{s} has strong Allee effects and the black curve is the case when \eqref{s} has weak Allee effects.}} \label{fig2:H}
   \end{figure}

\begin{figure}[ht]
\centering
   \includegraphics[scale =.65] {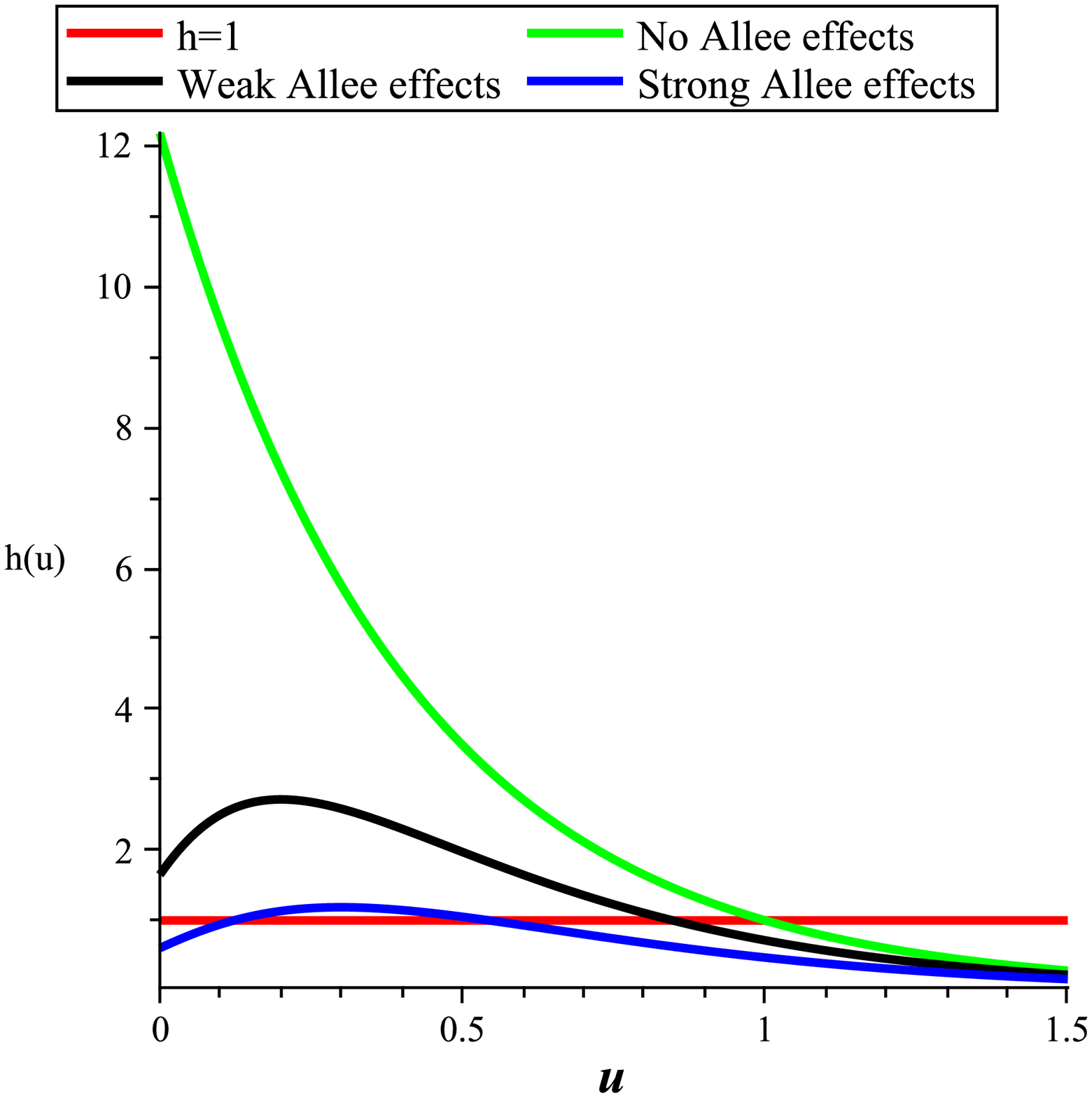}
  \caption{\modifyc{An illustration of strong Allee effects v.s. weak Allee effects of the per capita growth rate $h(u)=\frac{H(u)}{u}$:corresponding to its growth function $H(u)$ in Figure \ref{fig2:H}:  The red curve is $H=u$ or $h=1$; the green curve is the case when \eqref{s} has no Allee effects; the blue curve is the case when \eqref{s} has strong Allee effects and the black curve is the case when \eqref{s} has weak Allee effects.} }\label{fig2:h}
\end{figure}

For convenience, we define the following four conditions with respect to Model \eqref{s}:
\begin{itemize}
\item[]\textbf{A1:} There exists some $\epsilon>0$ such that $\frac{\partial h}{\partial u}>0$ for $u\in (0, \epsilon)$. 
\item[]\textbf{A2:} $ h(u)\geq 0$ and $\lim_{u\rightarrow\infty}h(u)=a<1$.
\item[]\textbf{A3:} There exists unique two numbers $0<A<K$, such that $h(A)=h(K)=1$ and $h'(A)>1$. In addition, $h(0)<1$.
\item[]\textbf{A4:} There exists a unique fixed point $K$, i.e., $h(K)=1$. In addition, $h(0)>1$.
\end{itemize}\modifyb{Condition \textbf{A1} indicates that the per capita growth rate has a positive correlation with population density for small populations. Condition \textbf{A2} indicates that the population is positive invariant (i.e., all future population is nonnegative for any nonegative initial population) and bounded. Condition \textbf{A3} indicates that the population has two steady states $A$ and $K$ where $A$, called \emph{Allee threshold}, is always unstable and $K$, called \emph{carry capacity}, can be locally stable under certain conditions. Condition \textbf{A4} indicates that species $u$ is persistent and may be locally stable at its carry capacity $K$}. Now we can define the following terms:
\begin{itemize}
\item Species $u$ has \emph{Allee effects} if its population model \eqref{s} satisfies Condition \textbf{A1, A2}.
\item Species $u$ has \emph{strong Allee effects} if its population model \eqref{s} satisfies Condition \textbf{A1, A2, A3} (see blue curves in Figure \ref{fig2:H}-\ref{fig2:h}).
\item Species $u$ has \emph{weak Allee effects} if its population model \eqref{s} satisfies Condition \textbf{A1, A2, A4} (see black curves in Figure \ref{fig2:H}-\ref{fig2:h}).
\end{itemize}
Our definition of \emph{Allee effects} is more general than Livadiotis and Elaydi (2012) since they only define \emph{strong Allee effects} as the following three conditions: i) $h'(u)>0$ for $u\in (0,\epsilon)$ for some $\epsilon>0$; ii) $h(0)<1$ and iii) There exists a unique $K>0$, such that $h(K)=1, h'(K)<0$. These two definitions can be considered the same. When species $u$ is subject to \emph{strong Allee effects}, i.e., its population model \eqref{s} satisfies Condition \textbf{A1, A2, A3}, its population dynamics can be summarized as follows: When initial population density is below $A$, the population goes extinct; When initial population density is above $A$, the population may be able to sustain; when $h(0)$ is large enough, then \eqref{s} may have \emph{the essential extinction} (i.e., for any positive initial condition, the population converges to 0 with probability 1) due to strong Allee effects (Schreiber 2003). See Figure \ref{fig2:H}-\ref{fig2:h} for \modifyc{an illustration of strong Allee effects (blue curves in both figures) versus weak Allee effects (black curves in both figures)}.\\


Outcomes of varied patterns of intra-specific competition (e.g., \emph{contest} competitions or \emph{scramble} competitions) have different consequences for population dynamics \eqref{s} when species experiences  \emph{Allee effects}. A well-known contest competition model introduced by Thomson (1993) for fisheries is called the Sigmoid Beverton-Holt model that can be described as
\bae\label{sA}
u_{t+1}&=&H(u_t)=u_th(u_t)=u_t\frac{r u^{\delta-1}_t}{1+u_t^\delta}
\eae where $r,\delta>0$. Model \eqref{sA} is a depensatory generalization of the Beverton-Holt stock-recruitment relationship used to develop a set of constraints designed to safeguard against overfishing (Thomson 1993; Myers \emph{et al} 1995; Myers and Mertz 1998; Stoner and Ray-Culp 2000; Myers 2001; Gascoigne and Lipcius 2004; Harry \emph{et al} 2012). One important feature of \eqref{sA} is that it exhibits the \emph{Allee effects} if $\delta>1$. The dynamical properties of \eqref{sA} have been proposed by Harry \emph{et al} (2012). Here we summarize the important dynamics of \eqref{sA} that will be used later for our two-species model as the following proposition:

\begin{proposition}[Dynamical of \eqref{sA}]\label{th1:sA}Let $r,\delta>0$ and \eqref{sA} has an initial condition $x_0 > 0$. If $\delta>1$, we define $r_{crit}=\delta (\delta-1)^{\frac{1}{\delta}-1}$. Then the following statements are true:
\begin{enumerate}
\item If $\delta \in (0, 1)$, then Model \eqref{sA} has two non-negative equilibria: an unstable equilibrium 0 and a stable positive equilibrium $K$.
\item If $\delta>1$ and $r < r_{crit}$, then the only equilibrium of Model \eqref{sA} is 0 and it is locally stable.
\item If  $\delta>1$ and $r > r_{crit}$, then Model \eqref{sA} has three non-negative equilibria: $0,\, A$ and $K$ such that $0 < A < \frac{r(\delta-1)}{\delta} < K$. The 0 equilibrium is stable with the basin of attraction $[0, A )$; $A$ is unstable (repellor), while $K$ is stable with the basin of attraction $(A ,\infty)$.
\end{enumerate}
\end{proposition}
\noindent\textbf{Notes:} Harry \emph{et al} (2012) gave similar results of Proposition \ref{th1:sA} as their Proposition 1but without rigorous proof. Here, we provide a detailed proof in the Appendix. Proposition \ref{th1:sA} implies that Model \eqref{sA} has \emph{Allee effects} for any $r, \delta>0$. Precisely, Model \eqref{sA} has \emph{weak Allee effects} if $\delta \in (0, 1)$ while Model \eqref{sA} has \emph{strong Allee effects} if $\delta>1$ and $r > r_{crit}$. In addition, we would like to point out that when $\delta=2$, Elaydi and Sacker (2009) have extended \eqref{sA} to the following generalized form:
\bae\label{sA-es}
u_{t+1}&=&u_t\frac{d u_t+e}{u_t^2+bu_t+c}
\eae where $b$ is called the \emph{shock recovery} parameter and $c$ can be considered as a carrying capacity.  The dynamics of \eqref{sA-es} is similar to \eqref{sA} when $\delta=2$. Define a single species population model as follows
\bae\label{sA-a}
u_{t+1}&=&H(u_t)=u_th(u_t)=u_t\frac{r u^{\delta-1}_t}{a+u_t^\delta}
\eae where $r$ is the maximum intrinsic growth rate and $a$ can be considered as a carrying capacity as the parameter $c$ in Model \eqref{sA-es}. Then the dynamics of Model \eqref{sA-a} can be obtained as the following corollary by directly applying the results from Proposition \ref{th1:sA} if we define $r_{crit}=a^{\delta}\delta (\delta-1)^{\frac{1}{\delta}-1}$:
\begin{corollary}[Dynamical of \eqref{sA-a}]\label{c1:sA-a}Let $r,a,\delta>0$ and \eqref{sA-a} has an initial condition $x_0 > 0$. If $\delta>1$, we define $r_{crit}^a=a^{1/\delta}\delta (\delta-1)^{\frac{1}{\delta}-1}.$ Then the results of Proposition \ref{th1:sA} can directly apply to Model \eqref{sA-a}. In particular, if $a,\delta>1$ and $r>r_{crit}^a$, then Model \eqref{sA-a} has two interior equilibria $A^{a}$ and $K^a$ such that $0<A<A^a<K^a<K$ where $A,K$ are defined in Proposition \ref{th1:sA}.
\end{corollary}

The scramble intra-specific competition models subject to \emph{Allee effects} can take a form of the Ricker's model, for example, 
\bae\label{sA-r}
u_{t+1}&=&u_th(u_t)=u_te^{r(1-u_t)-\frac{m}{1+bu_t}}
\eae where $r$ represent the intrinsic growth rate (also the scramble intra-specific competition coefficient); $m$ represents predation intensities and $b$ represents the product of the proportional to the handling time and the carrying capacity. For Model \eqref{sA-r}, we have
$$h(u)=e^{r(1-u)-\frac{m}{1+bu}} \mbox{ and } h'(u)=\frac{[bm-r(1+bu)^2]e^{r(1-u)-\frac{m}{1+bu}}}{(1+bu)^2}.$$
Thus, according to the definitions above, Model \eqref{sA-r} has \emph{Allee effects} if $mb>0$; Model \eqref{sA-r} has \emph{weak Allee effects} if $bm>r>m$; and Model \eqref{sA-r} has \emph{strong Allee effects} if $r<m<\frac{r(1+b)^2}{4b}, b>1$.
In addition, Model \eqref{sA-r} has more complicate dynamics (e.g., chaos, essential extinctions) than the contest competition model subject to Allee effects \eqref{sA}. The detailed dynamics of \eqref{sA-r} has been studied by Schreiber (2003). \\

Another example of scramble competition with Allee effects is a generalized Maynard Smith-Slatkin model that can be described as follows:
\bae\label{sA-mss}
u_{t+1}&=&H(u_t)=u_th(u_t)=u_t\frac{r u_t^{\delta-1}}{1+u_t^\delta+b u_t^{d}}
\eae where all parameters are positive; $r, \delta$ have the same biological meanings as in Model \eqref{sA} and $b, d$ have the same biological meaning as in Model \eqref{sA-es}. Model \eqref{sA-mss} can exhibit similar complicated dynamics as Model \eqref{sA-r} due to the generic feature of \emph{scramble competition}. Here we summarize the basic dynamics of \eqref{sA-mss} that will be used in the later sections as the following theorem:
\begin{theorem}[Dynamics of Model \eqref{sA-mss}]\label{th2:sA-mss} If $\delta>1$, we define $r_{crit}$ as follows $$r_{crit}=\left(\frac{\delta}{b(d-\delta)}\right)^{(1-\delta)/d}\left[\frac{d}{d-\delta}+\left(\frac{\delta}{b(d-\delta)}\right)^{\delta/d}\right].$$Depending on the values of $r, b, \delta, d$, the basic dynamics of \eqref{sA-mss} can be classified as the following cases:
\begin{enumerate}
\item If $\delta>d$, then the trajectory of Model \eqref{sA-mss} with initial condition $x_0>0$ converges to one of its equilibria. More specifically, we have the following two cases:
\begin{enumerate}
\item If $\delta<1$, then Model \eqref{sA-mss} has two non-negative equilibria 0 and $K$ where 0 is unstable and $K$ is globally stable.
\item If $\delta>1$, then $u=0$ is unstable and Model \eqref{sA-mss} can have 0 or two or even four positive equilibria. Moreover, Model \eqref{sA-mss} is globally stable at 0 when it has no positive interior equilibrium.
\end{enumerate}
\item If $\delta=d$, then Model \eqref{sA-mss} is reduced to Model \eqref{sA} where its dynamics has been stated in Proposition \ref{th1:sA}.
\item If $\delta<d$, then $u=0$ is always locally asymptotically stable. In addition, 
   \begin{enumerate}
   \item If $\delta<1$, then Model \eqref{sA-mss} has two non-negative equilibria 0 and $K$ where 0 is unstable.
   \item If $\delta>1$ and $r>r_{crit}$, then Model \eqref{sA-mss} has at least two distinct positive roots $A,K$ such that $0<A<u_c<K$.
     \end{enumerate}
\item If both $\delta$ and $d$ are positive integers, then Model \eqref{sA-mss} has at most two distinct positive interior equilibria.\end{enumerate}
\end{theorem}

\begin{figure}[ht]
\centering
   \includegraphics[scale =0.45] {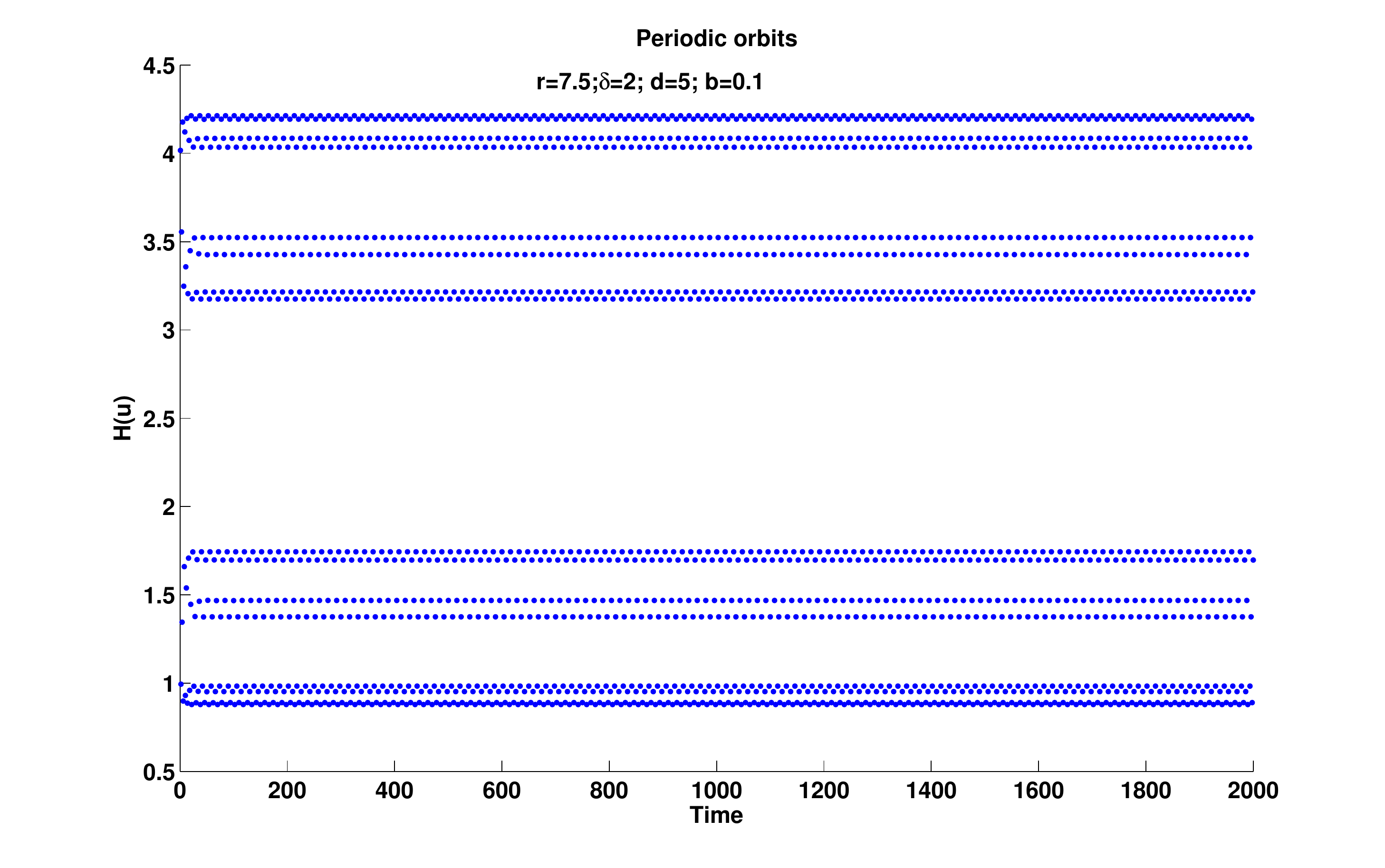} 
     \caption{\modifyc{Periodic time series of Model \eqref{sA-mss} when $\delta=2;d=5;b=0.1$ and $r=7.5$.} }\label{fig3:p} 
\end{figure}

\begin{figure}[ht]
\centering
   \includegraphics[scale =.45] {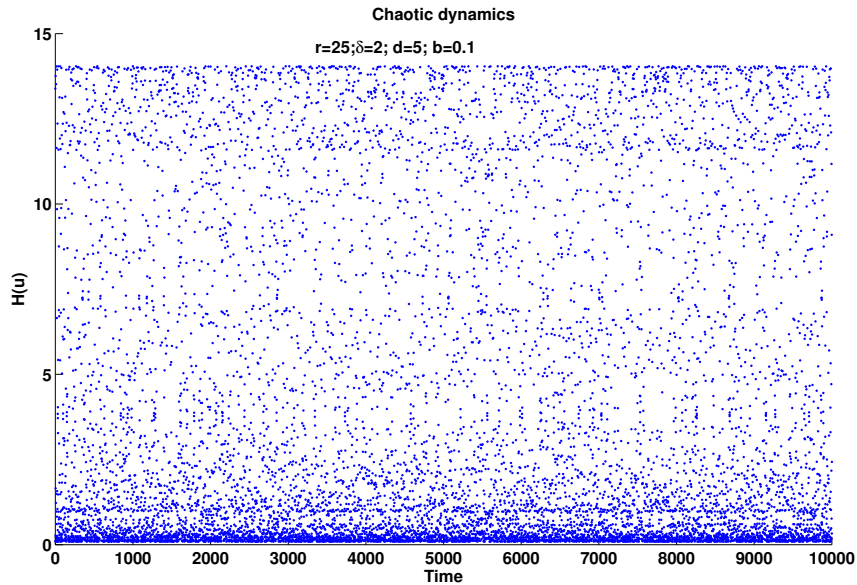}
  \caption{\modifyc{Chaotic time series of Model \eqref{sA-mss} when $\delta=2;d=5;b=0.1$ and $r= 25$.} }\label{fig3:c} 
\end{figure}

\noindent\textbf{Notes:} Theorem \ref{th2:sA-mss} \modifyb{implies that Model \eqref{sA-mss} can exhibit \emph{strong Allee effects} if $\delta>\max\{1,d\}$ or $1<\delta<d, r>r_{crit}.$ Model \eqref{sA} can exhibit \emph{strong Allee effects} only if $\delta>1$ and $r>r_{crit}.$} Even though Part 1 (b) and Part 3 (b) of Theorem \ref{th2:sA-mss} do not state that Model \eqref{sA-mss} has exactly two interior equilibria, numerical simulations suggest that Model \eqref{sA-mss} has at most two interior equilibria. In addition, if $1<\delta<d$ and $r>r_{crit}$, then \modifyr{Model \eqref{sA-mss} is a unimodal map that can have very complicated dynamics as the Ricker's map including \emph{the essential extinction} due to \emph{strong Allee effects}.} Numerical simulations suggest that increasing values of $r, b,\delta$ can destabilize the system while increasing value of $d$ can stabilize the system (see Figure \ref{fig3:c}-\ref{fig3:p} as examples).

\subsection{A generalized Beverton-Holt competition model with Allee effects}
\modifyb{A generalized Beverton-Holt two species competition model with Allee effects from single species models \eqref{sA} and \eqref{sA-mss} can be represented by the following two equations:} 
\bae\label{gx}
x_{t+1}&=&\frac{r_1 x_t^{\delta_1}}{1+x_t^{\delta_1}+b_1y_t^{\delta_3}}\\
\label{gy}
y_{t+1}&=&\frac{r_2 y_t^{\delta_2}}{1+y_t^{\delta_2}+b_2x_t^{\delta_4}}
\eae where $r_i, i=1,2$ describes the maximum intrinsic growth rate of each species; $b_i, i=1,2$ describes the relative inter-specific competition coefficient to its intra-specific competition; $\delta_i, i=1,2$ describes the nonlinearity of the intra-specific competition, while $\delta_i, i=3,4$ describes the nonlinearity of the inter-specific competition. Define Condition \textbf{H1}:
\begin{itemize}
\item[]\textbf{H1:} $\delta_i>1,\,r_i>r^i_{crit},\, i=1,2$ where $r^i_{crit}=\delta_i (\delta_i-1)^{\frac{1}{\delta_i}-1}, i=1,2$.
\end{itemize}
According to Proposition \ref{th1:sA}, \modifyb{each species of Model \eqref{gx}-\eqref{gy} suffers both contest intra-specific competition and strong Allee effects if Condition \textbf{H1} is satisfied}. Due to our interests, we study the population dynamics of \eqref{gx}-\eqref{gy} when Condition \textbf{H1} is satisfied for most cases.


\section{\modifyb{Dynamical properties of a generalized Beverton-Holt competition model with Allee effects}}
In this section, we explore the basic dynamical property of Model \eqref{gx}-\eqref{gy} in terms of its monotonicity, the extinction and coexistence of two competing species. 
If Condition \textbf{H1} is satisfied, Model \eqref{gx}-\eqref{gy} always has the following five boundary equilibria:
$$E_{0}=(0,0),\, E_{x_a0}=(A_1,0),\, E_{x_k0}=(K_1,0),\, E_{0y_a}=(0,A_2) \mbox{ and } E_{0y_k}=(0,K_2)$$ where $A_i, i=1,2$ is the \emph{Allee threshold} of species $i$ and $K_i, i=1,2$ is the carrying capacity of species $i$, thus, $0<A_i<K_i, i=1,2$. Simple calculations indicate that $E_0, E_{x_k0}$ and $E_{0y_k}$ are locally asymptotically stable, thus these three boundary equilibria are also attractors for Model \eqref{gx}-\eqref{gy}. This gives the following lemma:

\begin{lemma}[Basic properties]\label{l1:bp} \modifyr{With respect the map defined by Model \eqref{gx}-\eqref{gy} the positive cone $\mathbb R^2_+$ is forward invariant and the rectangular $[0, r_1]\times[0,r_2]$ is absorbing.}
If Condition \textbf{H1} is satisfied, then it has five boundary equilibria $E_0, E_{x_a0}, E_{x_k0}, E_{0y_a}$ and $E_{0y_k}$ where
$E_0, E_{x_k0}$ and $E_{0y_k}$ are locally asymptotically stable and $E_{x_a0},E_{y_a0}$ are unstable. In addition, let
$X_1=\{(x,0)\in\mathbb R^2_+\}$ and $X_2=\{(0,y)\in\mathbb R^2_+\}$, then both $X_i, i=1,2$ are positive invariant sets for Model \eqref{gx}-\eqref{gy}.
\end{lemma}
\noindent\textbf{Notes:} \modifyr{It is easy to check that the set $[0, r_1]\times[0,r_2]$ is forward invariant and all orbits in the positive cone $\mathbb R^2_+$ lie in this set after one time step. The stability of the boundary can be obtained by straight forward calculations, i.e.,  calculating the eigenvalues of the Jacobian matrix evaluated at the boundary equilibria.} Thus, we omit the details here. Also, we would like to point out that $E_{x_a0}$ and $E_{0y_a}$ are saddle nodes: $E_{x_a0}$ is a saddle node with the unstable manifold lying on the x-axis and the stable manifold tangent to the vertical eigenvector emanating from $E_{x_a0}$; similarly, $E_{0y_a}$ is a saddle node with the unstable manifold lying on the y-axis and the stable manifold tangent to the horizontal eigenvector emanating from $E_{0y_a}$.

\subsection{Competitive systems}
Monotonicity properties have been used successfully to study the dynamical behavior of solutions to differential equations (e.g., see [Smith 1995] for an overview) and difference equations (e.g., Selgrade and Ziehe 1987; Dancer and Hess 1991; Hess and Lazer 1991; Smith 1998; Wang and Jiang 2001; Clark \emph{et al} 2003). A competitive map is one type of monotone systems.\\

To continue our study, we first define a partial order on $\mathbb R^2$ so that the positive cone in this new partial order is the fourth quadrant. Specifically, for $u=(u_1,u_2), v=(v_1,v_2)\in\mathbb R^2_+$, we say that $u\leq_K v$ if $u_1\leq v_1$ and $v_2\leq u_2$. Two points $u, v\in\mathbb R^2_+$ are said to be \emph{related} (or $K$-ordered) if $u\leq_K v$ or $v\leq_K u$. A strict inequality between points can be defined as $u<_Kv$ if $u<_Kv$ and $u\neq v$. A stronger inequality can be defined as $u<<_Kv$ if $u_1<v_1$ and $v_2<u_2$. A map $T : Int\, \mathbb R^2_+\rightarrow Int \, \mathbb R^2_+$ is \emph{competitive} (or $K$-order preserving) if $u\leq_K v$ implies that $T(u)\leq_K T(v)$ for all $u,v\in Int\, \mathbb R^2_+$.  A map $T : Int \,\mathbb R^2_+\rightarrow Int\, \mathbb R^2_+$ is  \emph{strictly competitive} (\emph{strongly competitive}) if $u<_Kv$ implies that $T(u)<_KT(v)$ ($T(u)<<_KT(v)$) for all $u, v\in Int\,\mathbb R^2_+$. Clearly, being related is an invariant under iteration of a strongly competitive map. Let $J$ be the Jacobian matrix of a map $T$, then $J$ is $K$-positive ($K$-strongly positive) if its diagonal entries are nonnegative (positive) and its off-diagonal entries are nonpositive (negative).\\

\begin{theorem}[Competitive systems case 1]\label{th3:cs} If $\delta_1\delta_2\geq\delta_3\delta_4$, then Model \eqref{gx}-\eqref{gy} is strongly competitive and \modifyr{all orbits in the positive cone $\mathbb R^2_+$ converge to an equilibrium}. 
\end{theorem}
\noindent\textbf{Notes:} The proof of Theorem \ref{th3:cs} use the monotone properties of Model \eqref{gx}-\eqref{gy}, which is provided in the Appendix. This theorem indicates that Model \eqref{gx}-\eqref{gy} has \modifyr{equilibrium dynamics} (e.g., no nontrivial periodic orbits) whenever $\delta_1\delta_2\geq\delta_3\delta_4$. Recall that $\delta_i, i=1,2$ \modifyr{measures the intensity} of the intra-specific competition of species, and $\delta_i, i=3,4$ represents  the inter-specific competition. Then we can \modifyr{define $\delta_1\delta_2, \delta_3\delta_4$ as intensities of the inter- and intra-specific competition, respectively. Then the result is that the dynamics equilibration when the intensity of intra-specific competition intensity exceeds that of inter-specific competition.}\\


From the proof of Theorem \ref{th3:cs}, we can see that Model \eqref{gx}-\eqref{gy} is \emph{strictly competitive} in $\mathbb R^2_+$ and is \emph{strongly competitive} in $Int \,\mathbb R^2_+$. For the case that $\delta_1\delta_2<\delta_3\delta_4$, we have the following theorem:\\

\begin{theorem}[Competitive systems case 2]\label{th5:cs2} Assume that $\delta_1\delta_2<\delta_3\delta_4$. If 
$$r_1<\left(\frac{\delta_1\delta_2}{b_2[\delta_3\delta_4-\delta_1\delta_2]}\right)^{1/\delta_4}\mbox{ or }r_2<\left(\frac{\delta_1\delta_2}{b_1[\delta_3\delta_4-\delta_1\delta_2]}\right)^{1/\delta_3},$$then every orbit of Model \eqref{gx}-\eqref{gy} converges to an equilibrium.
\end{theorem}

\noindent\textbf{Notes:} The proof of Theorem \ref{th5:cs2} use the monotone properties of Model \eqref{gx}-\eqref{gy} as well as the properties of non-invertible maps. The detailed of proof has been given in the Appendix.
In the case that $\delta_1\delta_2<\delta_3\delta_4$, Theorem \ref{th5:cs2} indicates that the dynamics of \eqref{gx}-\eqref{gy} can still \modifyr{have equilibrium dynamics} as the case when $\delta_1\delta_2\geq\delta_3\delta_4$ (Theorem \ref{th3:cs}) if the maximum intrinsic growth rate of one species is not too large. \modifyb{However, when this additional condition is not satisfied, Model \eqref{gx}-\eqref{gy} can have \modifyr{non-equilibrium attractors} as suggested by results of Teresc\'ak (1996) and Smith (1998). 
For example, if we let $\delta_1=\delta_2=\delta, \delta_3=\delta_4=d, r_1=r_2, b_1=b_2$, then the dynamics of \eqref{gx}-\eqref{gy} on the invariant manifold $\Omega_{x=y}=\{(x,y)\in\mathbb R^2_+: x=y\}$ can be reduced to Model \eqref{sA-mss} which can have complicated dynamics. \modifyc{See Figure \ref{fig3:p}-\ref{fig3:c} as examples for different dynamics on $\Omega_{x=y}$ when $\delta=2;d=5;b=0.1$ and $r=7.5$ (Figure \ref{fig3:p}) and $r= 25$ (Figure \ref{fig3:c})}}. \\


\subsection{Basins of attractions of \modifyr{extinction equilibria}}
Assume that Condition \textbf{H1} is satisfied for Model \eqref{gx}-\eqref{gy}. We can define the following sets:
$$O_0=[0,A_1]\times[0,A_2],\,\, O_{ex}=\{(x,y)\in \mathbb R^2_+: x<A_1\} \mbox{ and } O_{ey}=\{(x,y)\in \mathbb R^2_+: y<A_2\}$$
and Condition \textbf{H2}-\textbf{H3} as follows:
\begin{itemize}
\item[]\textbf{H2:} $\delta_1>1$ and  $r_1>r^{a_1}_{crit}$ where $r^{a_1}_{crit}=(a_1)^{1/\delta_1}\delta_1 (\delta_1-1)^{\frac{1}{\delta_1}-1}$ and $a_1=1+b_1(A_2)^{\delta_3}$.
\item[]\textbf{H3:} $\delta_2>1$ and  $r_2>r^{a_2}_{crit}$ where $r^{a_2}_{crit}=(a_2)^{1/\delta_2}\delta_2 (\delta_2-1)^{\frac{1}{\delta_2}-1}$ and $a_2=1+b_2(A_1)^{\delta_4}$.
\end{itemize}
\modifyc{Condition \textbf{H2} or \textbf{H3} indicates that species suffers from strong Allee effects.} Assume that Model \eqref{sA-a} satisfies Condition \textbf{H2}, then according to Corollary \ref{c1:sA-a}, it has two interior equilibria $A_1^{a_1}$ and $K_1^{a_1}$ where $A_1^{a_1}$ is the Allee threshold and $K_1^{a_1}$ is the carrying capacity. Similarly, if Model \eqref{sA-a} satisfies Condition \textbf{H3}, then it has two interior equilibria $A_1^{a_2}$ and $K_1^{a_2}$ where $A_2^{a_2}$ is the Allee threshold and $K_2^{a_2}$ is the carrying capacity. Let 
$$O_{x}=\{(x,y)\in \mathbb R^2_+: x>A_1^{a_1},\, y<A_2\} \mbox{ and } O_{y}=\{(x,y)\in \mathbb R^2_+: x<A_2,\, y>A_2^{a_2}\},$$ then we can show the following theorem:\\
\begin{theorem}[Basins of attractions of boundary attractors I]\label{th4:bs} Assume that Model \eqref{gx}-\eqref{gy} satisfies Condition \textbf{H1}, then
\begin{itemize}
\item For any initial value in $O_{ex}$, we have\,\,
$\lim_{t\rightarrow\infty}x_t=0.$
\item For any initial value in $O_{ey}$, we have\,\,
$\lim_{t\rightarrow\infty}y_t=0.$
\item For any initial value in $O_0$, we have\,\,
$\lim_{t\rightarrow\infty}(x_t, y_t)=E_0.$
\end{itemize}
If Model \eqref{sA-a} satisfies Condition \textbf{H2}, then for any initial value in $O_{x}$, we have
$$\lim_{t\rightarrow\infty}(x_t,y_t)=E_{x_{k_1}0}.$$
If Model \eqref{sA-a} satisfies Condition \textbf{H3}, then for any initial value in $O_{y}$, we have
$$\lim_{t\rightarrow\infty}(x_t,y_t)=E_{0y_{k_2}}.$$

\end{theorem}

\noindent\textbf{Notes:} Theorem \ref{th4:bs} (see the proof in the Appendix) provides an approximation of the boundary attractor's basins of attractions for Model \eqref{gx}-\eqref{gy} when each species suffers from contest competition and strong Allee effects. \modifyb{In the absence of other species, each species' extinction region is $[0, A_i]$ while in the presence of other species, the extinction region becomes larger due to the inter-specific competition, thus, competition cannot save species from extinction at their low densities. In order to maintain species x's population, its initial condition should be large enough and the other species population should be small enough, e.g., the initial condition should be larger than $A_1^{a_1}>A_1$ while species y's initial condition should be less than $A_2$.} In addition, Theorem \ref{th4:bs} implies that the basins of attractions of the boundary attractors $E_{x_{k_1}0}$ and $E_{0y_{k_2}}$ are unbounded (see Figure \ref{fig5:bn}, \ref{fig5:b2i} and \ref{fig5:b1i}).\\

Define 
$$\begin{array}{lcl}
O_{ex}^l&=&\bigcup_{k>0}\{(x,kx)\in Int\,\mathbb R^2_+: x>\left(\frac{r_1-A_1}{b_1k^{\delta_3}}\right)^{\frac{1}{\delta_3-\delta_1}}\}\\
O_{ey}^l&=&\bigcup_{k>0}\{(ky,y)\in Int\,\mathbb R^2_+: y>\left(\frac{r_2-A_2}{b_2k^{\delta_4}}\right)^{\frac{1}{\delta_4-\delta_2}}\}\end{array}$$ and let
$O_0^l=O_{ex}^l\cap O_{ey}^l.$ Then we have the following theorem:
\begin{theorem}[Basins of attractions of boundary attractors II]\label{th5:bsE0} Assume that Model \eqref{gx}-\eqref{gy} satisfies Condition \textbf{H1}.
\begin{itemize}
\item If $\delta_3>\delta_1$, then for any initial value in $O_{ex}\cup O_{ex}^l$, we have\,\,
$\lim_{t\rightarrow\infty}x_t=0.$
\item If $\delta_4>\delta_2$, then for any initial value in $O_{ey}\cup O_{ey}^l$, we have\,\,
$\lim_{t\rightarrow\infty}y_t=0.$
\item If $\delta_3>\delta_1$ and $\delta_4>\delta_2$, then for any initial value in $O_0\cup O_0^l$, we have\,\,
$\lim_{t\rightarrow\infty}(x_t, y_t)=E_0.$
\end{itemize}
\end{theorem}
\noindent\textbf{Notes:} The detailed proof of Theorem \ref{th5:bsE0} is provided in the Appendix. \modifyb{This theorem indicates that if a species' relative competition degree is less than 1 (i.e., $\delta_3>\delta_1$ or $\delta_4>\delta_2$), then its extinction region (i.e., the initial conditions that lead to the extinction of two species) consists of two distinct components, of which one is stated in Theorem \ref{th4:bs} and the other one is the area that when its competitor's population is large enough. This phenomenon is is a typical property of non-invertible maps which is caused by the large degree of the inter-specific competition, i.e., $\delta_3>\delta_1$ or $\delta_4>\delta_2$.} Theorem \ref{th5:bsE0} also indicates that the basins of attractions of $E_0$ may be unbounded. In fact, the basins of attractions of $E_0$ is a connected and unbounded region when there is no interior equilibrium (see Figure \ref{fig5:n} and \ref{fig5:bn}). While the basins of attractions of $E_0$ consists of two parts when there is an interior equilibrium: one part is a bounded region including the neighborhood of $E_0$ while the other part is a connected unbounded region including very large initial values of both species (see Figure \ref{fig5:b2i} and \ref{fig5:b1i}).

\subsection{\modifyr{Persistent} equilibria}
\begin{figure}[ht]
\centering
 \includegraphics[scale =0.45] {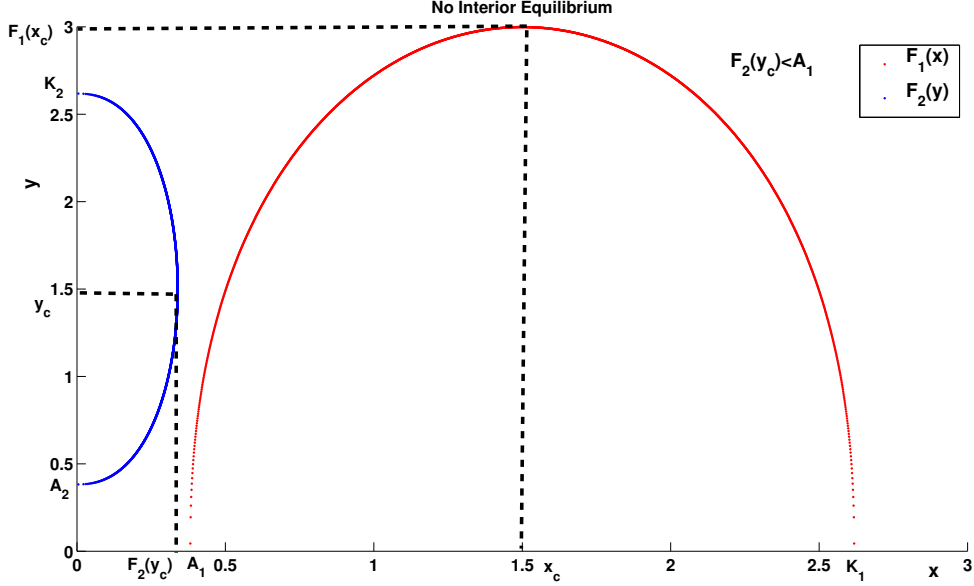}
 \caption{Schematic nullclines when Model \eqref{gx}-\eqref{gy} has no interior equilibrium of Case 1, i.e., $F_2(y_c)<A_1$: \modifyc{The red curve is $F_1(x)=\left(\frac{r_1x^{\delta_1-1}-x^{\delta_1}-1}{b_1}\right)^{1/\delta_3}$ where its maximum $F_1(x_c)$ occurs at $x_c=\frac{r_1(\delta_1-1)}{\delta_1}$; While the blue curve is $F_2(y)=\left(\frac{r_2y^{\delta_2-1}-y^{\delta_2}-1}{b_2}\right)^{1/\delta_4}$ where its maximum $F_2(y_c)$ occurs at $y_c=\frac{r_2(\delta_2-1)}{\delta_2}$}. } \label{fig4:no}
\end{figure}

\begin{figure}[ht]
\centering
\includegraphics[scale =.45] {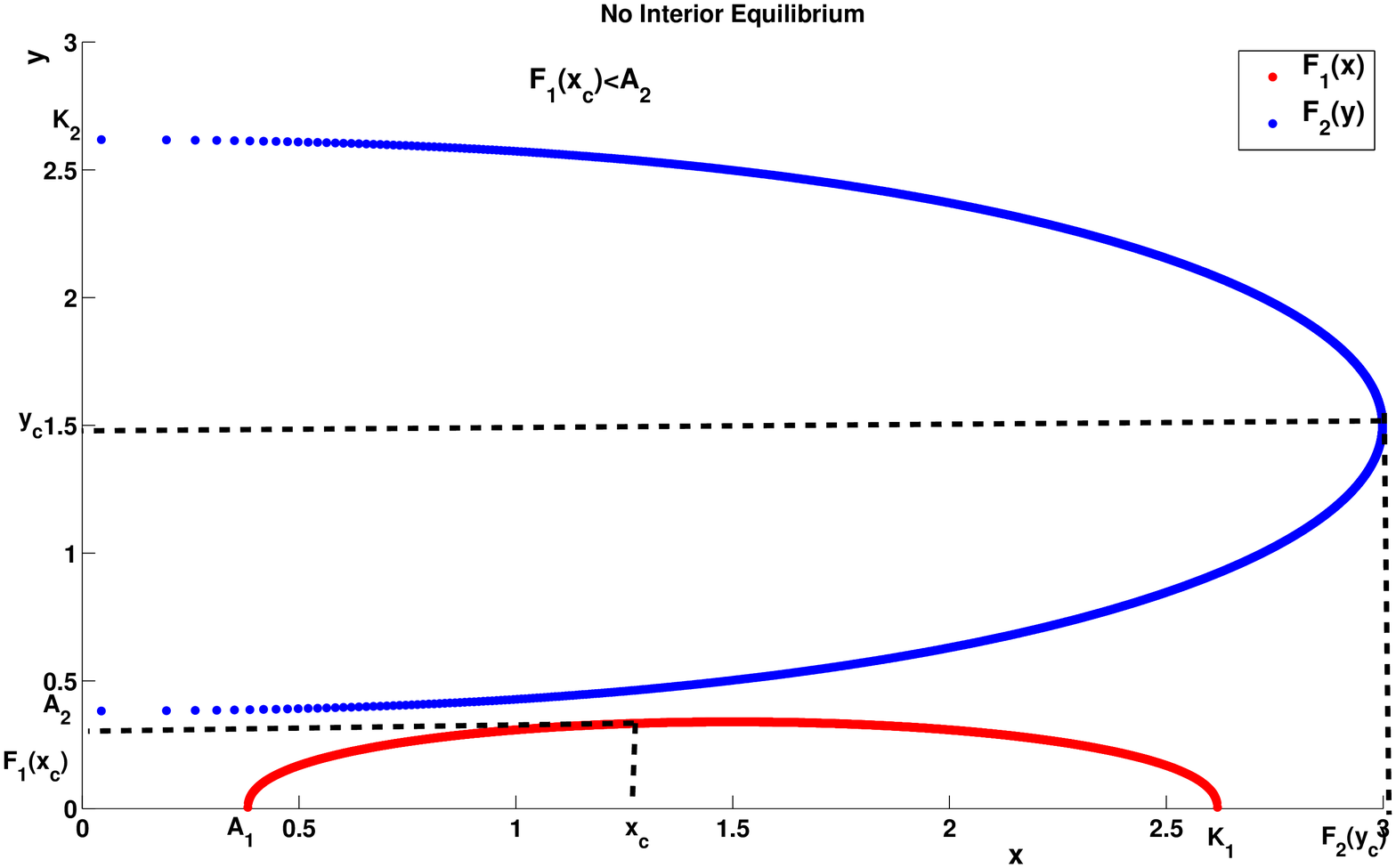}
 \caption{Schematic nullclines when Model \eqref{gx}-\eqref{gy} has no interior equilibrium of Case 2, i.e., $F_1(x_c)<A_2$: \modifyc{The red curve is $F_1(x)=\left(\frac{r_1x^{\delta_1-1}-x^{\delta_1}-1}{b_1}\right)^{1/\delta_3}$ where its maximum $F_1(x_c)$ occurs at $x_c=\frac{r_1(\delta_1-1)}{\delta_1}$; While the blue curve is $F_2(y)=\left(\frac{r_2y^{\delta_2-1}-y^{\delta_2}-1}{b_2}\right)^{1/\delta_4}$ where its maximum $F_2(y_c)$ occurs at $y_c=\frac{r_2(\delta_2-1)}{\delta_2}$}. } \label{fig4:no2}
\end{figure}

\begin{figure}[ht]
\centering
 \includegraphics[scale =0.45] {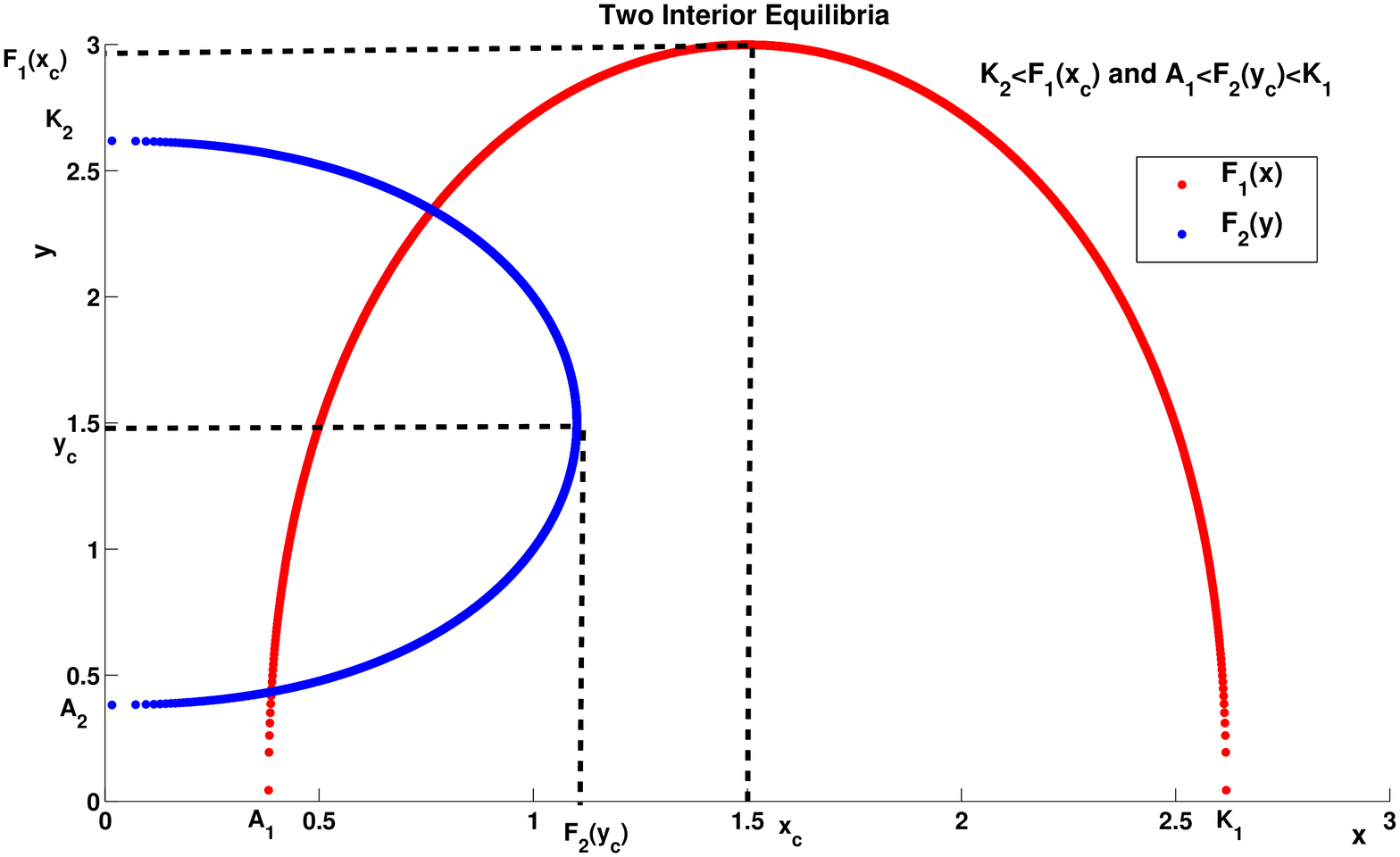}
  \caption{Schematic nullclines when Model \eqref{gx}-\eqref{gy} has two interior equilibria of Case 1, i.e., $K_2<F_1(x_c), \,A_1<F_2(y_c)<K_1$: \modifyc{The red curve is $F_1(x)=\left(\frac{r_1x^{\delta_1-1}-x^{\delta_1}-1}{b_1}\right)^{1/\delta_3}$ where its maximum $F_1(x_c)$ occurs at $x_c=\frac{r_1(\delta_1-1)}{\delta_1}$; While the blue curve is $F_2(y)=\left(\frac{r_2y^{\delta_2-1}-y^{\delta_2}-1}{b_2}\right)^{1/\delta_4}$ where its maximum $F_2(y_c)$ occurs at $y_c=\frac{r_2(\delta_2-1)}{\delta_2}$}. } \label{fig4:2i}
\end{figure}

\begin{figure}[ht]
\centering
 \includegraphics[scale =0.45] {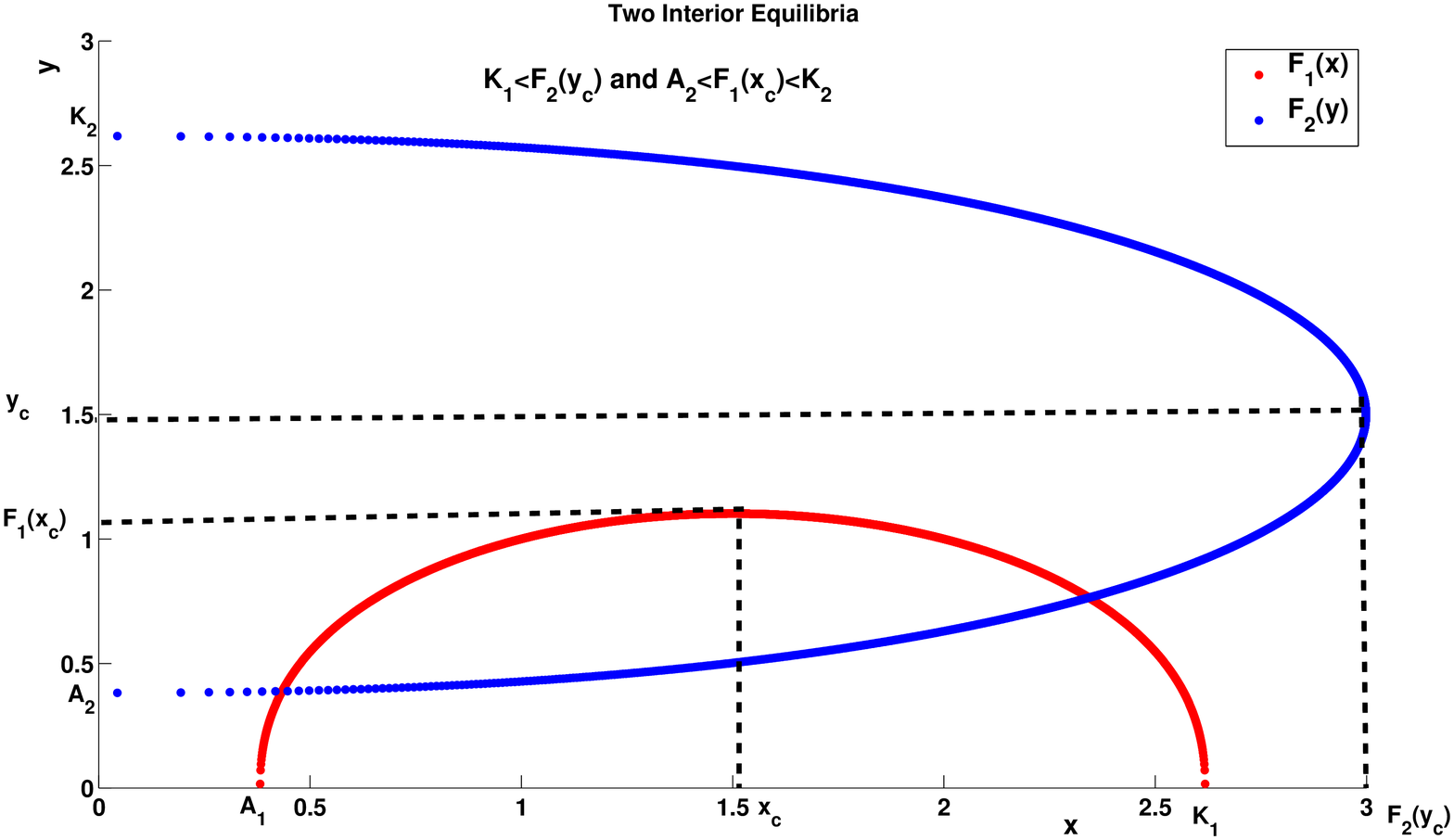}
  \caption{Schematic nullclines when Model \eqref{gx}-\eqref{gy} has two interior equilibria of Case 2, i.e., $A_2<F_1(x_c)<K_2, \,K_1<F_2(y_c)$: \modifyc{The red curve is $F_1(x)=\left(\frac{r_1x^{\delta_1-1}-x^{\delta_1}-1}{b_1}\right)^{1/\delta_3}$ where its maximum $F_1(x_c)$ occurs at $x_c=\frac{r_1(\delta_1-1)}{\delta_1}$; While the blue curve is $F_2(y)=\left(\frac{r_2y^{\delta_2-1}-y^{\delta_2}-1}{b_2}\right)^{1/\delta_4}$ where its maximum $F_2(y_c)$ occurs at $y_c=\frac{r_2(\delta_2-1)}{\delta_2}$}. } \label{fig4:2i2}
\end{figure}

\begin{figure}[ht]
\centering
 \includegraphics[scale =0.45] {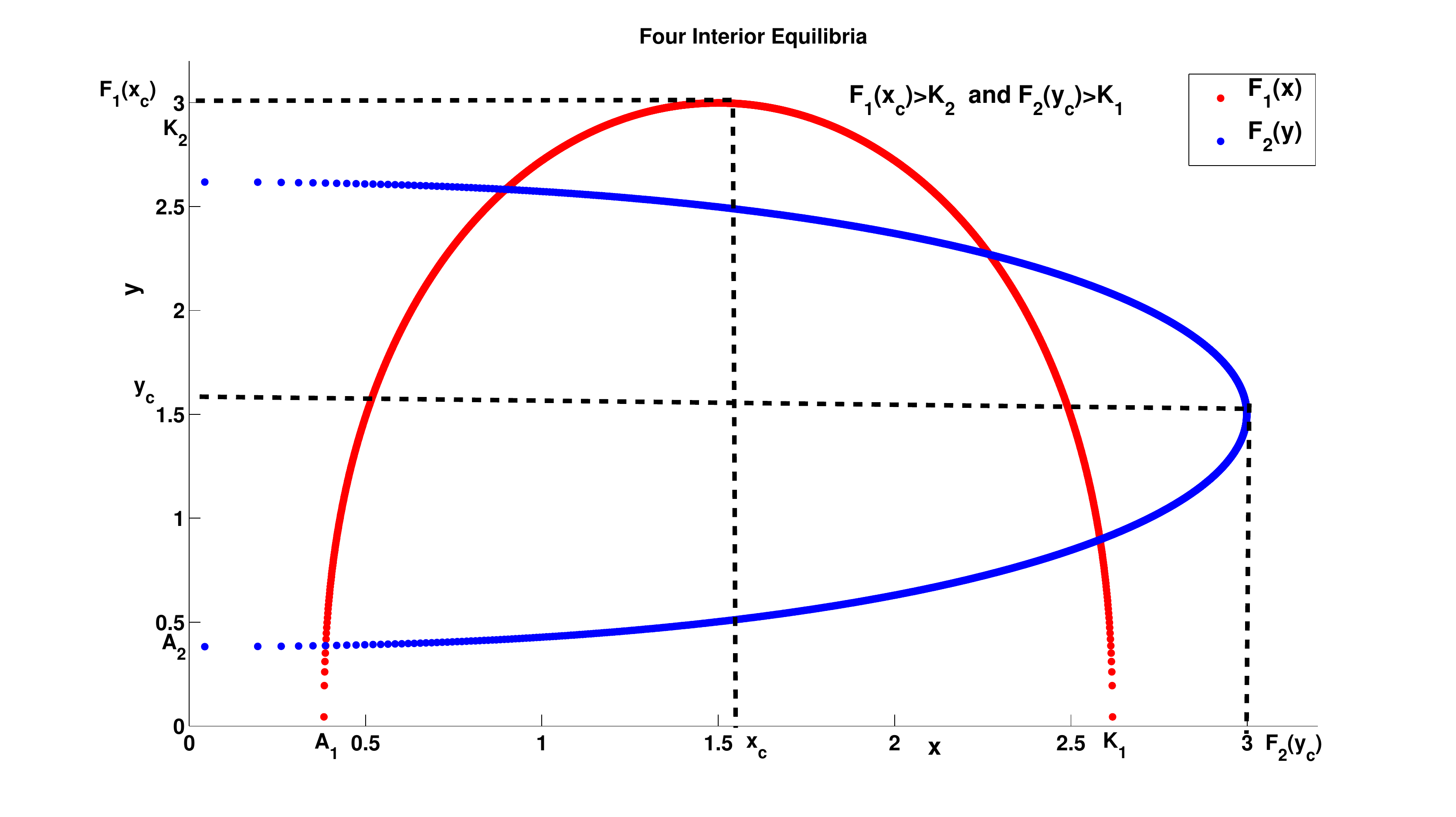}
  \caption{Schematic nullclines when Model \eqref{gx}-\eqref{gy} has four interior equilibria, i.e., $F_1(x_c)>K_2$ and $F_2(y_c)>K_1$: \modifyc{The red curve is $F_1(x)=\left(\frac{r_1x^{\delta_1-1}-x^{\delta_1}-1}{b_1}\right)^{1/\delta_3}$ where its maximum $F_1(x_c)$ occurs at $x_c=\frac{r_1(\delta_1-1)}{\delta_1}$; While the blue curve is $F_2(y)=\left(\frac{r_2y^{\delta_2-1}-y^{\delta_2}-1}{b_2}\right)^{1/\delta_4}$ where its maximum $F_2(y_c)$ occurs at $y_c=\frac{r_2(\delta_2-1)}{\delta_2}$}. }\label{fig4:4i}
\end{figure}

Assume that Model \eqref{gx}-\eqref{gy} satisfies Condition \textbf{H1}, then according to Theorem \ref{th4:bs}, it always has three boundary attractors. This suggests that the coexistence of two species is possible only if there exists an interior equilibrium. We are going to prove the following theorem regarding the number of interior equilibria.
\begin{theorem}[Interior equilibria]\label{th6:ie}Assume that Model \eqref{gx}-\eqref{gy} satisfies Condition \textbf{H1}. Let 
$$F_1(x)=\left(\frac{r_1x^{\delta_1-1}-x^{\delta_1}-1}{b_1}\right)^{1/\delta_3}\mbox{ and }F_2(y)=\left(\frac{r_2y^{\delta_2-1}-y^{\delta_2}-1}{b_2}\right)^{1/\delta_4},$$ then we have the following statements:
\begin{itemize}
\item If $F_1(x_c)<A_2$ or $F_2(y_c)<A_1$, then Model \eqref{gx}-\eqref{gy} has no interior equilibrium.
\item If $A_2<F_1(x_c)<K_2, \,K_1<F_2(y_c)$ or $A_1<F_2(y_c)<K_1,\,K_2<F_1(x_c)$, then Model \eqref{gx}-\eqref{gy} has two interior equilibria.
\item If $F_1(x_c)>K_2$ and $F_2(y_c)>K_1$, then Model \eqref{gx}-\eqref{gy} has four interior equilibria.
\end{itemize}where $x_c=\frac{r_1(\delta_1-1)}{\delta_1}$ and $y_c=\frac{r_2(\delta_2-1)}{\delta_2}$.
\end{theorem}
\noindent\textbf{Notes:} \modifyr{See Figure \ref{fig4:no}-\ref{fig4:4i} regarding the generic nullclines of Model \eqref{gx}-\eqref{gy}. Theorem \ref{th6:ie} (its proof stated in the Appendix) gives us two scenarios when Model \eqref{gx}-\eqref{gy} has interior equilibria where coexistence occurs only if Model \eqref{gx}-\eqref{gy} has four interior equilibria (see Figure \ref{fig5:1i} and Figure \ref{fig5:b1i}). Sufficient conditions guaranteeing four interior equilibria is that $F_1(x_c)>K_2$ and $F_2(y_c)>K_1$, i.e. $$x_{c}<\frac{r_1 x_c^{\delta_1}}{1+x_c^{\delta_1}+b_1K_2^{\delta_3}}\mbox{ and } y_{c}<\frac{r_2 y_c^{\delta_2}}{1+y_c^{\delta_2}+b_2K_1^{\delta_4}}$$ which indicates that any initial condition in $[x_c, K_1]\times[y_c,K_2]$ lead to a locally asymptotically stable interior equilibrium according to the monotonicity of the system. Biologically, this implies that both species are able to persist under proper initial conditions if $F_1(x_c)>K_2$ and $F_2(y_c)>K_1$ holds.}
In the next subsection, we explore the local stability of a symmetric system.



\subsection{A symmetric competition model with Allee effects}
In order to obtain more insights on the coexistence of two competing species, we focus on a symmetric case of \eqref{gx}-\eqref{gy} where $r_1=r_2=r, b_1=b_2=b, \delta_1=\delta_2=\delta, \delta_3=\delta_4=d$ which gives the following model
\bae\label{sgx}
x_{t+1}&=&\frac{r x_t^\delta}{1+x_t^\delta+by_t^{d}}\\
\label{sgy}
y_{t+1}&=&\frac{r y_t^\delta}{1+y_t^\delta+b x_t^{d}}. 
\eae
Understanding the dynamics of \eqref{sgx}-\eqref{sgy} also can help us obtain a better understanding of the dynamics of asymmetric cases due to the robust property of our system \eqref{gx}-\eqref{gy}.

 It is easy to verify that \eqref{sgx}-\eqref{sgy} is reduced to the 1-D model described by \eqref{sA-mss} if $x_0=y_0$, i.e., the manifold defined by $\Omega_{y=x}=\{(x,y)\in\mathbb R^2_+:x=y\}$ is invariant.
According to Theorem \ref{th2:sA-mss}, Model \eqref{sgx}-\eqref{sgy} has at least two symmetric equilibria $(A, A)$ and $(K,K)$ such that $0<A<\left(\frac{r(d-\delta)}{\delta}\right)^{1/d}<K$ if 
\bae\label{stability}d>\delta>1\mbox{ and }r>r_{crit}=\left(\frac{\delta}{b(d-\delta)}\right)^{(1-\delta)/d}\left[\frac{d}{d-\delta}+\left(\frac{\delta}{b(d-\delta)}\right)^{\delta/d}\right].\eae
Now we have the following theorem regarding the stability of $(K,K)$.
\begin{theorem}[Stability of the symmetric equilibrium]\label{th7:sie} Assume that Model \eqref{sgx}-\eqref{sgy} satisfies the inequalities \eqref{stability}. Then \eqref{sgx}-\eqref{sgy} has at least two symmetric equilibria $(A, A)$ and $(K,K)$ such that $0<A<\left(\frac{r(d-\delta)}{\delta}\right)^{1/d}<K$. Moreover, the symmetric interior equilibrium $(x^*,x^*)$ is locally asymptotically stable if 
\bae\label{st}\Big\vert\frac{\delta+b (\delta-d) (x^*)^d}{1+(x^*)^{\delta}+b(x^*)^d}\Big\vert<1 \mbox{ and } \frac{\delta+b (\delta+d)(x^*)^d}{1+(x^*)^{\delta}+b(x^*)^d} <1.\eae
While the symmetric interior equilibrium $(x,x)$ is unstable if 
\bae\label{us}\Big\vert\frac{\delta+b (\delta-d) (x^*)^d}{1+(x^*)^{\delta}+b(x^*)^d}\Big\vert>1 \mbox{ or } \frac{\delta+b (\delta+d)(x^*)^d}{1+(x^*)^{\delta}+b(x^*)^d} >1.\eae
\end{theorem}

\noindent\textbf{Notes:} \modifyb{Theorem \ref{th7:sie} indicates that it is impossible for the interior equilibrium of the symmetric model \eqref{sgx}-\eqref{sgy} to have Neimak-Sacker bifurcation from the invariant $\Omega_{x=y}$ since the eigenvalue that determines the points moving towards or away from $\Omega_{x=y}$ is always positive (see the detailed proof in the Appendix). This suggests that \eqref{sgx}-\eqref{sgy} may have relatively simple dynamics}.

\section{Application to an example}
In this section, we let $\delta_1=\delta_2=2$, then Model \eqref{gx}-\eqref{gy} can be rewritten as follows:
\bae\label{x}
x_{t+1}&=&\frac{r_1 x_t^2}{1+x_t^2+b_1 y_t^{\delta_3}}\\
\label{y}
y_{t+1}&=&\frac{r_2 y_t^2}{1+y_t^2+b_2 x_t^{\delta_4}}.
\eae In the absence of species y, i.e., $y_0=0$, we have \eqref{x}-\eqref{y} is reduced to a single species model \eqref{sA} when $r=r_1$ and $\delta=2$, i.e., 
\bae\label{sAx}
x_{t+1}&=&\frac{r_1 x_t^2}{1+x_t^2}.
\eae \modifyb{Model \eqref{sAx} is a well-known fisheries model when fish is subject to contest intra-specific competition and strong Allee effects.} The dynamics of \eqref{sAx} is very simple: 1. If $r_1<2$, then \eqref{sAx} converges to 0 for any $x_0>0$. 2. If $r_1>2$, then \eqref{sAx} has two interior equilibria $A_1$ and $K_1$ where
$$A_1=\frac{r_1-\sqrt{r_1^2-4}}{2} \mbox{ and } K_1=\frac{r_1+\sqrt{r_1^2-4}}{2}$$ such that \eqref{sAx} converges to 0 if $x_0<A_1$ while  \eqref{sAx} converges to $K_1$ if $x_0>A_1$.

Similarly, in the absence of species x, i.e., $x_0=0$, \eqref{x}-\eqref{y} is reduced to a single species model \eqref{sA} for the case that $r=r_1$ and $\delta=2$, i.e., 
\bae\label{sAy}
y_{t+1}&=&\frac{r_2 y_t^2}{1+y_t^2}
\eae where it also has two interior equilibria $A_2$ and $K_2$ where
$$A_2=\frac{r_2-\sqrt{r_2^2-4}}{2} \mbox{ and } K_2=\frac{r_2+\sqrt{r_2^2-4}}{2}.$$
Due to our interests, we focus on the case that $r_i>2, i=1,2$ for Model \eqref{x}-\eqref{y}. Then apply Theorem \ref{th3:cs}, Theorem \ref{th4:bs}, Theorem \ref{th5:cs2}, Theorem \ref{th5:bsE0} and Theorem \ref{th6:ie}, we obtain the following corollary regarding the dynamics of \eqref{x}-\eqref{y}:

\begin{corollary}[Dynamics of \eqref{x}-\eqref{y}]\label{c2:ds2} Assume that $r_i>2$. Define 
$$r^{a_1}_{crit}=\sqrt{1+b_1\left(\frac{r_2-\sqrt{r_2^2-4}}{2}\right)^{1/\delta_3}} \mbox{ and }r^{a_2}_{crit}=\sqrt{1+b_2\left(\frac{r_1-\sqrt{r_1^2-4}}{2}\right)^{1/\delta_4}}.$$
 Then the following statements are true:
\begin{enumerate}
\item If $\delta_3\delta_4\leq4$, then every orbit of Model \eqref{x}-\eqref{y} with any initial condition in $\mathbb R^2_+$ converges to one of its equilibria. 
\item If $\delta_3\delta_4>4$ and 
$$r_1<\left(\frac{4}{b_2[\delta_3\delta_4-4]}\right)^{1/\delta_4}\mbox{ or }r_2<\left(\frac{4}{b_1[\delta_3\delta_4-4]}\right)^{1/\delta_3},$$then every orbit of Model \eqref{x}-\eqref{y} with any initial condition in $\mathbb R^2_+$ converges to one of its equilibria. 
\item Model \eqref{x}-\eqref{y} always has three boundary attractors $E_0, E_{x_{k_1}0}, E_{0y_{k_2}}$ where
$$E_0=(0,0),\, E_{x_{k_1}0}=(\frac{r_1+\sqrt{r_1^2-4}}{2},0) \mbox{ and } E_{0y_{k_2}}=(0,\frac{r_2+\sqrt{r_2^2-4}}{2}).$$ Moreover, 
\begin{itemize}
\item For any initial value in $O_{ex}$, we have\,\,
$\lim_{t\rightarrow\infty}x_t=0.$ While if $\delta_3>2$, then for any initial value in $O_{ex}\cup O_{ex}^l$, we have\,\,
$\lim_{t\rightarrow\infty}x_t=0.$
\item For any initial value taken in $O_{ey}$, we have\,\,
$\lim_{t\rightarrow\infty}y_t=0.$ While if f $\delta_4>2$, then for any initial value in $O_{ey}\cup O_{ey}^l$, we have\,\,
$\lim_{t\rightarrow\infty}y_t=0.$
\item  For any initial value taken in $O_0$, we have\,\,
$\lim_{t\rightarrow\infty}(x_t, y_t)=E_0.$ While if $\delta_3>2$ and $\delta_4>2$, then for any initial value in $O_0\cup O_0^l$, we have\,\,
$\lim_{t\rightarrow\infty}(x_t, y_t)=E_0.$
\item If $r_1>r^{a_1}_{crit}$, then for any initial value in $O_{x}$, we have
$$\lim_{t\rightarrow\infty}(x_t,y_t)=E_{x_{k_1}0}.$$ While if $r_2>r^{a_2}_{crit}$, then for any initial value in $O_{y}$, we have
$$\lim_{t\rightarrow\infty}(x_t,y_t)=E_{0y_{k_2}}.$$
\end{itemize}
\item Model \eqref{x}-\eqref{y} has no interior equilibrium if 
$$\left(\frac{r_1^2/4-1}{b_1}\right)^{1/\delta_3}<\frac{r_2-\sqrt{r_2^2-4}}{2}\mbox{ or }\left(\frac{r_2^2/4-1}{b_2}\right)^{1/\delta_4}<\frac{r_1-\sqrt{r_1^2-4}}{2}.$$Model \eqref{x}-\eqref{y} has two interior equilibria if 
$$\frac{r_2-\sqrt{r_2^2-4}}{2}<\left(\frac{r_1^2/4-1}{b_1}\right)^{1/\delta_3}<\frac{r_2+\sqrt{r_2^2-4}}{2}\mbox{ and }\left(\frac{r_2^2/4-1}{b_2}\right)^{1/\delta_4}>\frac{r_1+\sqrt{r_1^2-4}}{2}$$ or
$$\frac{r_1-\sqrt{r_1^2-4}}{2}<\left(\frac{r_2^2/4-1}{b_2}\right)^{1/\delta_4}<\frac{r_1+\sqrt{r_1^2-4}}{2}\mbox{ and }\left(\frac{r_1^2/4-1}{b_1}\right)^{1/\delta_3}>\frac{r_2+\sqrt{r_2^2-4}}{2}.$$Model \eqref{gx}-\eqref{gy} has four interior equilibria if 
$$\left(\frac{r_1^2/4-1}{b_1}\right)^{1/\delta_3}>\frac{r_2+\sqrt{r_2^2-4}}{2}\mbox{ and }\left(\frac{r_2^2/4-1}{b_2}\right)^{1/\delta_4}>\frac{r_1+\sqrt{r_1^2-4}}{2}.$$
\end{enumerate}
\end{corollary}
The corollary above does not give us information on the dynamic patterns of the system. Due to the complexity of the model, we use numerical simulations to investigate their dynamical patterns. According to simulations, we summarize its dynamical patterns as follows:

\begin{itemize}
\item \textbf{No interior equilibrium:} When Model  \eqref{x}-\eqref{y} has no interior equilibrium (see Figure \ref{fig5:n}),  it has exactly three boundary attractors, i.e., (0,0), $(x_2,0)$ and $(0,y_2)$, and no interior attractor. The basins of attractions of these three attractors are unbounded (shown in Figure \ref{fig5:bn}) where the extinction area (i.e., basins of attractions of $E_0=(0,0)$) is an open and unbounded region in red. \modifyb{This implies that the inter-specific competition make both two competing species prone to extinction.}
\begin{figure}[ht]
\centering
   \includegraphics[scale =0.45] {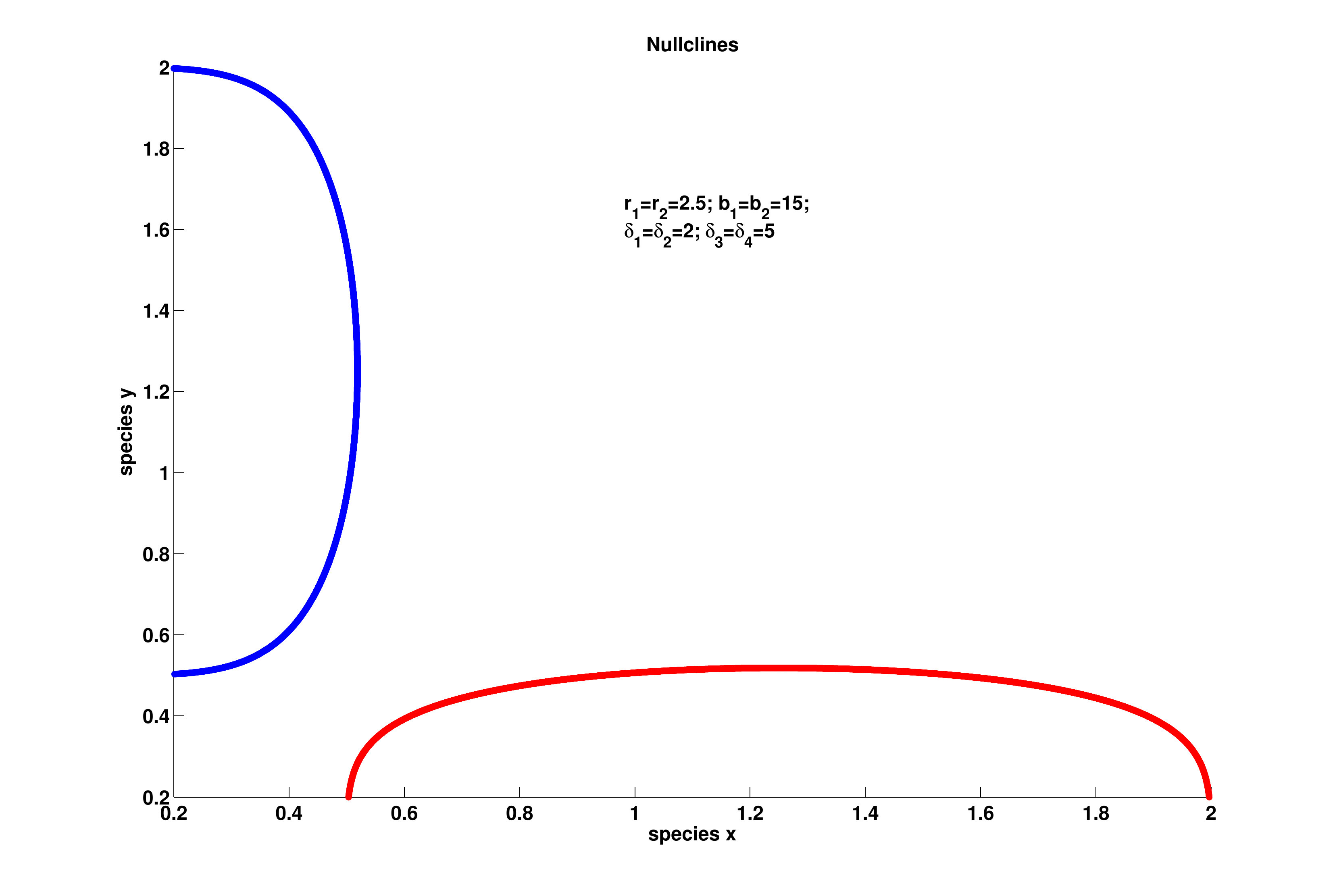}
   \caption{\modifyc{Nucllines of \eqref{x}-\eqref{y} when $r_1=r_2=2.5$, $b_1=b_2=15$ and $\delta_3=\delta_4=5$}.}  \label{fig5:n}
 \end{figure}
 \begin{figure}[ht]
\centering
   \includegraphics[scale =.45] {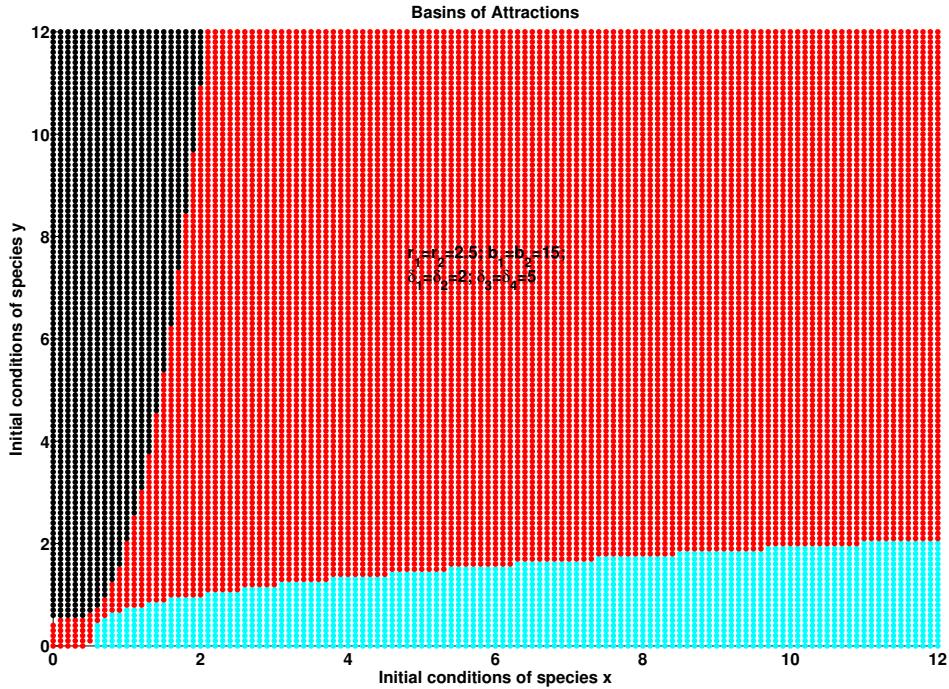}
   \caption{\modifyc{Basins of attractions of \eqref{x}-\eqref{y} when $r_1=r_2=2.5$, $b_1=b_2=15$ and $\delta_3=\delta_4=5$:  The black region is the basins of attractions of $(0,y_2)$; The cyan region is the basins of attractions of $(x_2,0)$; The blue region is the basins of attractions of the interior attractor which is an interior equilibrium and the red region is the basins of attractions of $(0,0)$.}}
    \label{fig5:bn}
\end{figure}
\item \textbf{Two interior equilibrium:} When Model  \eqref{x}-\eqref{y} has two interior equilibrium (see Figure \ref{fig5:2i}), it has exactly three boundary attractors, i.e., (0,0), $(x_2,0)$ and $(0,y_2)$, and no interior attractor. The basins of attractions of these three boundary attractors are unbounded (shown in Figure \ref{fig5:b2i}): The basins of attractions of $(x_2,0)$ and $(0,y_2)$ are connected and unbounded while the basins of attractions of $E_0$ consists of a bounded-connected and an unbounded-connected region.
\begin{figure}[ht]
\centering
   \includegraphics[scale =0.45] {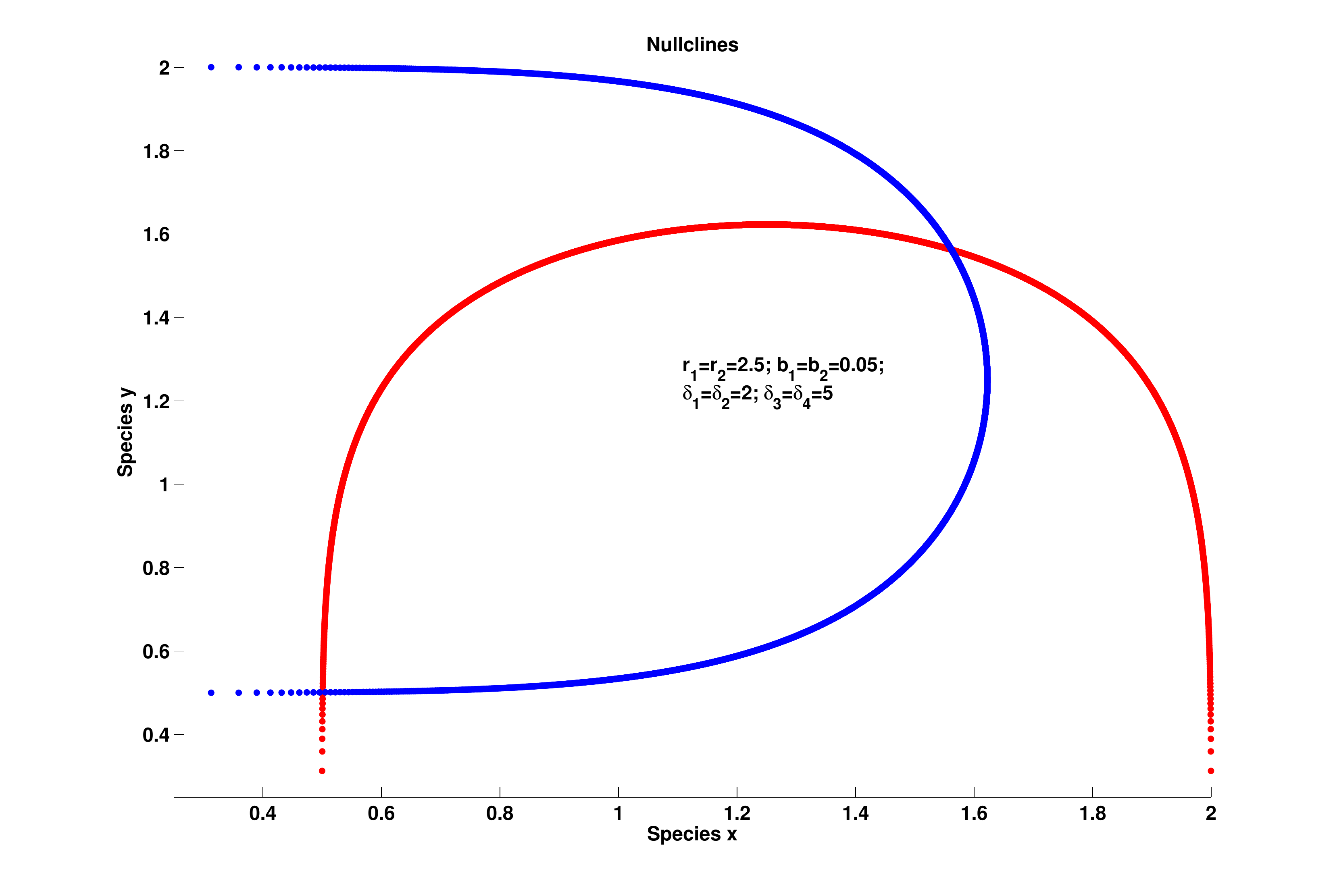} 
   \caption{\modifyc{Nucllines of \eqref{x}-\eqref{y} when $r_1=r_2=2.5$, $b_1=b_2=.05$ and $\delta_3=\delta_4=5$}.}\label{fig5:2i}
   \end{figure}
    \begin{figure}[ht]
\centering
   \includegraphics[scale =.45] {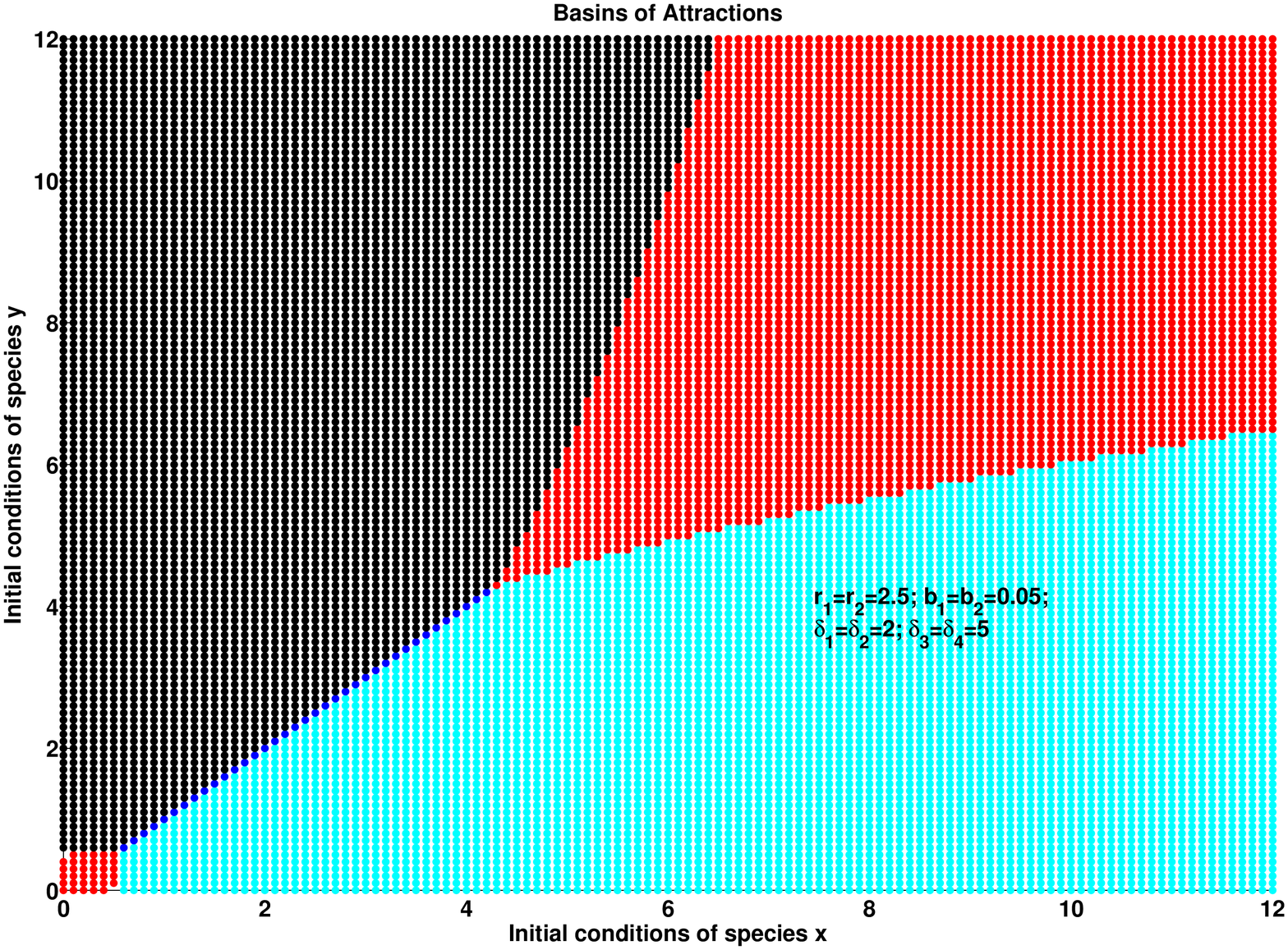} 
  \caption{\modifyc{Basins of attractions of \eqref{x}-\eqref{y} when $r_1=r_2=2.5$, $b_1=b_2=.05$ and $\delta_3=\delta_4=5$:  The black region is the basins of attractions of $(0,y_2)$; The cyan region is the basins of attractions of $(x_2,0)$; The blue region is the basins of attractions of the interior attractor which is an interior equilibrium and the red region is the basins of attractions of $(0,0)$.}}  \label{fig5:b2i}
\end{figure}
\item \textbf{Four interior equilibrium:} This is the only case when Model  \eqref{x}-\eqref{y} has an interior attractor (i.e., the coexistence of two species). When Model  \eqref{x}-\eqref{y} has four interior equilibrium (see Figure \ref{fig5:1i}), it has four attractors of which three are boundary attractors: (0,0), $(x_2,0)$ and $(0,y_2)$, and one is an interior attractor $(x^*,y^*)$ whose basins of attractions is a bounded and connected region. The basins of attractions of these four attractors are shown in Figure \ref{fig5:b1i} where the features of basins of attractions of boundary attractors are the same as the two interior equilibria case. 

\begin{figure}[ht]
\centering
   \includegraphics[scale =0.45] {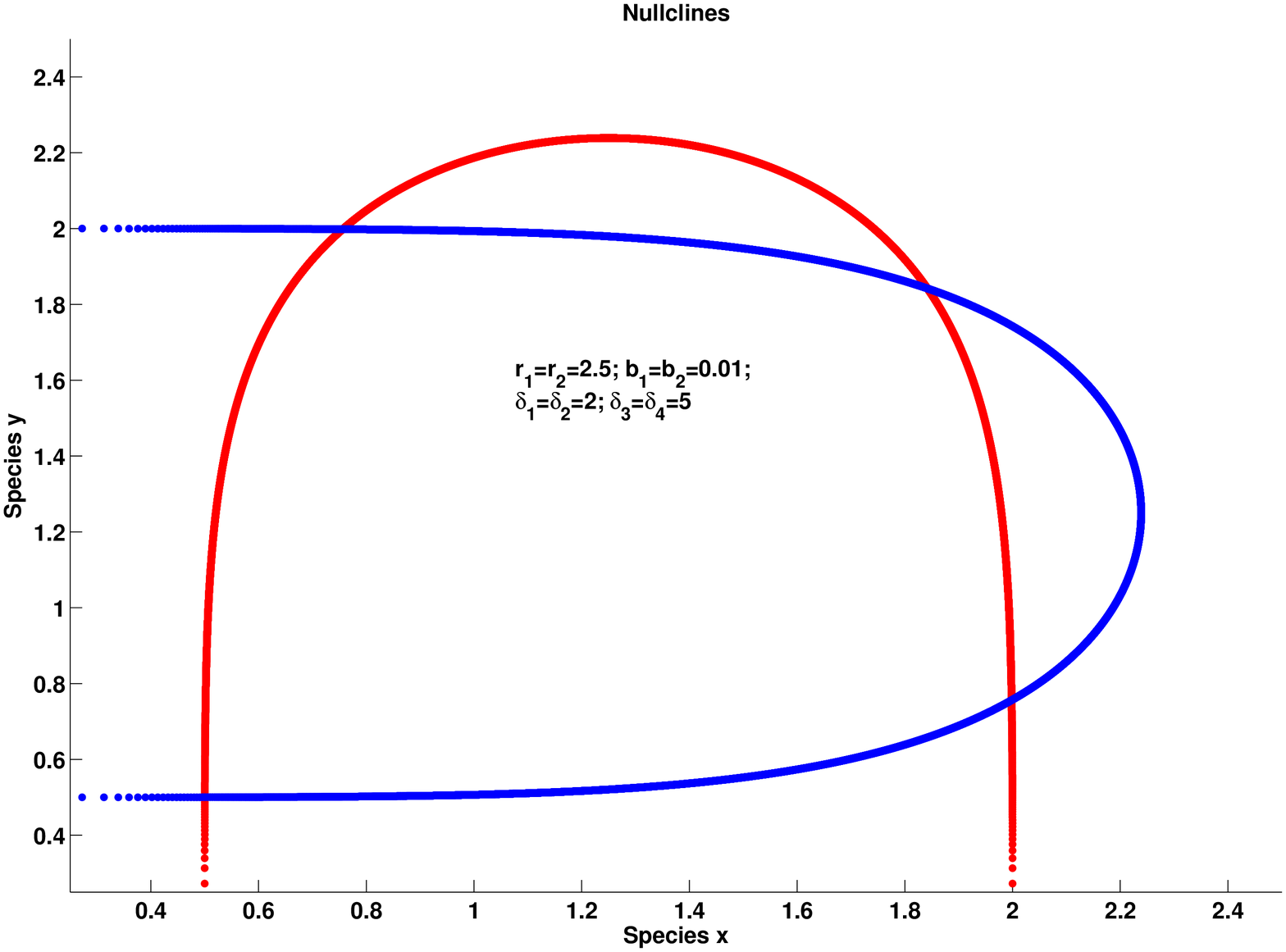}
    \caption{\modifyc{Nucllines of \eqref{x}-\eqref{y} when $r_1=r_2=2.5$, $b_1=b_2=.01$ and $\delta_3=\delta_4=5$.}}\label{fig5:1i}
\end{figure} 
   \begin{figure}[ht]
\centering
   \includegraphics[scale =.45] {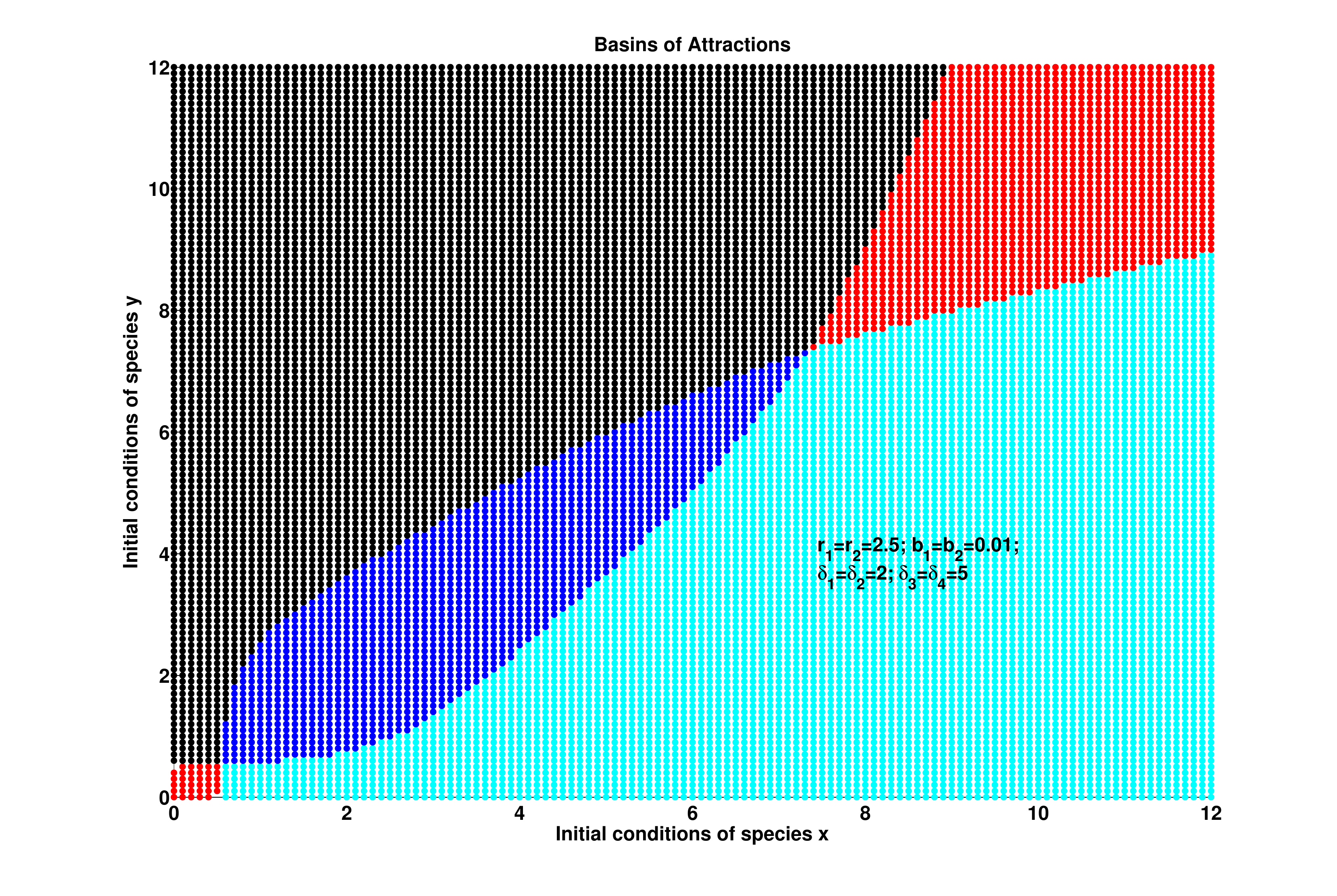} 
  \caption{\modifyc{Basins of attractions of \eqref{x}-\eqref{y} when $r_1=r_2=2.5$, $b_1=b_2=.01$ and $\delta_3=\delta_4=5$:  The black region is the basins of attractions of $(0,y_2)$; The cyan region is the basins of attractions of $(x_2,0)$; The blue region is the basins of attractions of the interior attractor which is an interior equilibrium and the red region is the basins of attractions of $(0,0)$.}} \label{fig5:b1i}
\end{figure}
\end{itemize} \modifyb{The observations above for Model \eqref{x}-\eqref{y} should be able to apply to a general model \eqref{gx}-\eqref{gy} due to the robustness of the model. In addition, we would like to point out that Model \eqref{gx}-\eqref{gy} can have coexistence of two species only if it has four interior equilibria.}

\section{Discussion}
\modifyc{Competition is an important ecological process which, in the short run, can cause a reduction in the number of species living within an area, by preventing very similar species from co-occurring; while in the long term, is likely to increase species diversity, by acting as a force for specialization and divergence. Competition can vary between two extreme forms: contest versus scramble. Allee effect is an another important ecological process that leads to a positive correlation of the per capita growth rate of a species and its low population density. The interplay of Allee effects and intra-specific competition can result in a critical threshold, called Allee threshold, below which species goes to extinction; while above which species may persist. This is referred to strong Allee effects. When there is no such threshold, it is referred as to weak Allee effects. Many species in nature are subject to both Allee effects and competition. The combinations of Allee effects and different type of competition can generate varied dynamical outcomes, especially for the coexistence of species.}\\

{ Including Allee effects in competition models can have huge potential impacts on diversity and community structure (Chesson and Ellner 1989; Hopf \emph{et al} 1993). Allee effects impose a cost of rarity (Mechanism R). When communities of organisms are subject to a cost of rarity imposed by Allee effects, they should be comprised of distinct species with differences in resource use that are greater than expected by chance (Hopf \emph{et al} 1993). Kang and Yakubu (2011) and Kang (2013) studied the population dynamics of a two competing species model where each species suffers from Allee effects and scramble competition: Kang and Yakubu (2011) showed that weak Allee effects may promote the permanence of two competing species at their low densities. The biological explanation for this is that weak Allee effects decrease the fitness of resident species so that the other species is able to invade at its low densities. Kang (2013) showed that strong Allee effects may save two competing species from essential extinctions at their high densities. The biological explanation for this is that scramble competition can bring the current high population density to a lower population density but is above the Allee threshold in the next season with the consequence that both competing species are able to persist. To extend the study by Kang and Yakubu (2011) and Kang (2013), this article studies the population dynamics of a generalized Beverton-Holt two species competition model where each species has strong Allee effects to explore the establishment and the extinction of two competing species subject to contest competition and Allee effects.}\\


\subsection{Two dimensional strong Allee effects of symmetric competition models}

Adopting the definition of Allee effects for single species, we are able to define Allee effects of a symmetric two species competition model that can be represented as 
\bae\label{sx}
x_{t+1}&=&F(x_t, y_t)=x_t f(x_t,y_t)\\
\label{sy}
y_{t+1}&=&F(y_t, x_t)=y_t f(y_t,x_t).
\eae where $\frac{\partial f(x,y)}{\partial y}<0$. Then we say two dimensional system \eqref{sx}-\eqref{sy} has \emph{two dimensional strong Allee effects} if both $f(x,0)$ and $f(x,x)$ satisfy Condition \textbf{A1, A2, A3}. Model \eqref{sgx}-\eqref{sgx} subject to inequalities \eqref{stability} provides an example of two competing species suffering from \emph{contest competition} and \emph{two dimensional strong Allee effects}. The ecological examples when two competing species suffering from \emph{scramble competition} and \emph{two dimensional strong Allee effects} are the following two examples:\\

\noindent\textbf{Model I:} Scramble competition with strong Allee effects due to mating limitations
\bae\label{sxs1}
x_{t+1}&=&x_t e^{r(1-x_t)-ay_t}\frac{b x_t}{1+bx_t}\\
\label{xys1}
y_{t+1}&=&y_t e^{r(1-y_t)-ax_t}\frac{b y_t}{1+by_t}
\eae where $r$ represent the intrinsic growth rates as well as intra-specific competition coefficient; $b$ represent the product of an individual's searching efficiency and the carrying capacities and $a$ represents the inter-specific competition coefficient. If $$\frac{e^{r-1}}{r+a}>\frac{1}{r+a}+\frac{1}{b},$$ then \textbf{Model I} has \emph{two dimensional strong Allee effects}.\\

\noindent\textbf{Model II:} Scramble competition with strong Allee effects due to predation saturations

\bae\label{sxs2}
x_{t+1}&=&x_t e^{r(1-x_t)-ay_t-\frac{m}{1+bx_t}}\\
\label{xys2}
y_{t+1}&=&y_t e^{r(1-y_t)-ax_t-\frac{m}{1+by_t}}
\eae  where $r$ represents the intrinsic growth rates as well as intra-specific competition coefficient; $m$ represents predation intensities; $b$ represent the product of the proportional to the handling time and the carrying capacities and $a$ represents the inter-specific competition coefficient. If $$r<m<\frac{(br+r+a)^2}{4b(a+r)} \mbox{ and } b>1+a/r,$$ then \textbf{Model I} has \emph{two dimensional strong Allee effects}. \\

\modifyc{Kang and Yakubu (2011) and Kang (2013) have studied the detailed dynamics of \textbf{Model I \& II} introduced above. By comparing our results to their study on a two competing species model where each species is subject to scramble competition and strong Allee effects (i.e., \textbf{Model I \& II}), we are able to obtain some generic dynamical features of a two-species competition model subject to strong Allee effects and the different dynamical outcomes due to different types of competition introduced in the models as follows:\\\vspace{5pt} }

\noindent\modifyc{\textbf{Generic dynamic features of scramble v.s. contest competition models with strong Allee effects:}
\begin{enumerate}
\item \modifyb{Both contest and scramble competition models with strong Allee effects can have either tri-stability (i.e., three boundary attractors $E_0, E_{x_{k_1}0}, E_{0y_{k_1}}$) or four attractors (i.e., three boundary attractors plus an interior attractor) depending on certain conditions.} 
\item The coexistence of two competing species is possible only if the system has at least four interior equilibria. 
\item Both models can exhibit \emph{two dimensional strong Allee effects} in the symmetric case.
\item \modifyb{The basins of attractions of $E_0$ are unbounded, which suggests that two competing species are not able to coexist if their population densities are too low or too high.}
 \end{enumerate}}

\noindent\modifyc{\textbf{Different dynamical outcomes of scramble v.s. contest competition models with strong Allee effects:} 
\begin{enumerate}
\item  \modifyb{Our generalized Beverton-Holt two species (representing contest competition) model} always has three boundary attractors $E_0, E_{x_{k_1}0}, E_{0y_{k_1}}$. While \modifyb{two-species Ricker models (representing scramble competition)} subject to strong Allee effects can have either one (i.e., extinction of two species) or three boundary attractors where the extinction is caused by strong Allee effects (i.e., each species has \emph{essential extinction} at its single state). 
\item Scramble competition models with strong Allee effects can promote the coexistence of two competing species at their high densities under certain condition. While this is not the case for contest competition models with strong Allee effects.
\item Scramble competition models can have complicated structures of attractors (e.g., chaos, \emph{essential extinctions}) while the contest competition model may have relative simple dynamics (e.g., from performed simulations, we only observed locally asymptotically stable equilibriua as attractors).
\item The basins of attractions of attractors that single species persists are bounded and may consist of different distinct parts for scramble competition models (see black and cyan areas in Figure 3\&4 in Kang 2013) while the basins of attractions of attractors that single species persists are unbounded and connected for contest competition models (see black and cyan areas in Figure \ref{fig5:bn}, \ref{fig5:b2i} and \ref{fig5:b1i}). 
\end{enumerate}}
\modifyb{The discussion above regarding contest v.s. scramble competition models with strong Allee effects are based on our study on the generalized Beverton-Holt competition model \eqref{gx}-\eqref{gy} and Ricker-type competition models with Allee effects studied in Kang and Yakubu (2011) and Kang (2013). The comparison may provide us useful insights of general two species competition models where each species has \emph{strong Allee effects} in its single state. However, in order to obtain more general results, we need to classify the different types of intra- and inter-specific competition rigorously and study how their combinations with Allee effects lead to different dynamical outcomes. For example, the following example can be considered as a combination of contest inter-specific competition and scramble intra-specific competition subject to Allee effects:
 $$\begin{array}{lcl}
 x_{t+1}&=&\frac{x_t e^{r_1(1-x_t)}I^x(x_t)}{1+c_1y_t}\\
  y_{t+1}&=&\frac{y_t e^{r_2(1-y_t)}I^y(y_t)}{1+c_2x_t}\end{array}$$where $I^i,i=x,y$ represents positive density dependence that can lead to Allee effects and $c_i,i=1,2$ are relatively inter-specific competition coefficients. This could be our future research project.}

\subsection{Conclusion}
\modifyc{In this article, we propose a generalized Beverton-Holt competition model to investigate how the interactions of Allee effects and contest competition affect the species's establishment and extinction. More precisely, we would like to answer the following questions addressed in the introduction: 1. What are the population dynamical outcomes of a two competing species when each of them suffers from both strong Allee effects and contest competition? 2. How does contest competition combined with strong Allee effects affect species' persistence and extinction by comparing to models subject to strong Allee effects and scramble competition? 3. What are the generic dynamical features of population models subject to strong Allee effects and competition? Here, we summarize our main results from the study as follows:
\begin{enumerate}
\item Theorem \ref{th3:cs} and Theorem \ref{th5:cs2} reveal that a two species contest competition model with Allee effects has equilibrium dynamics, under the conditions that \modifyr{the intensity of intra-specific competition intensity either exceeds that of inter-specific competition or is less than that of inter-specific competition but the maximum intrinsic growth rate of one species is not too large.}
 This partially answers the first question.
\item Theorem \ref{th4:bs} and Theorem \ref{th5:bsE0} provide us an approximations of basins of attractions of boundary attractors (i.e., initial conditions that lead to the extinction of one or both species). \modifyb{This indicates that the inter-specific competition can make species prone to extinction if it suffers from both strong Allee effects and contest competition, and the situation is worsen if species' relative competition degree is less than 1}. In addition, Theorem \ref{th6:ie}, \ref{th7:sie} combined with numerical simulations suggest that decreasing the values of the inter-specific competition coefficient and degree may promote the coexistence of two species. These results may answer the second question.
\item The third question can be partially answered by comparing with dynamics of two species scramble competition models subject to strong Allee effects in Kang (2013) to our Theorem \ref{th6:ie} combined with numerical simulations. The comparison suggests that (i) Both contest and scramble competition models can have either three boundary attractors or four attractor with only one interior attractor, however, scramble competition model may have only one attractor of the extinction of both species due to strong Allee effects; (ii) Scramble competition models represented by Ricker-type models may have more complicated structure of interior attractor than contest competition ones; and (iii) Scramble competition may be more likely to promote the coexistence of two species at low or high densities under certain conditions. 
\end{enumerate}
Our study combined with the results from Kang and Yakubo (2011) as well as Kang (2013) provide us useful insights on how Allee effects interacting with different forms of competition have an impact on the persistence and extinction of species. One direct application of these results is that we may be able to save rare species from \emph{essential extinction} due to strong Allee effects by introducing another similar competing species under proper conditions. Thus, our study may help us understand the diversity of nature and make better policy for conservation program, especially for rare species that suffers from Allee effects. 
In the future, it will be interesting to investigate the following questions: 1. What are population dynamics of two competing species when each of them or both have different form of competition in their different life stages as well as Allee effects? 2. What are necessary conditions for a general two species population model to generate two dimensional strong Allee effects? 3. How can we predict more precise initial conditions that lead to the extinction and the coexistence of two species?}

\section*{Acknowledgement}
The author would like to thank the useful discussions with Dr. Peter Chesson, Dr. Elaydi Saber, Dr. Lidicker William, Dr. Livadiotis George and Dr. Kwessi Eddy for the early draft of the manuscript. Y.K.'s research is partially supported by NSF-DMS (1313312), Simons Collaboration Grants for Mathematicians (208902) and the research scholarship from College of Letters and Sciences (CLS) at Arizona State University.
\section*{\modifyc{Appendix}}
\appendix
\subsection*{Proof of Proposition \ref{th1:sA}}
\begin{proof}First, it is easy to see that Model \eqref{sA} is positively invariant and bounded in $\mathbb R_+$. Notice that 
$$H(0)=0, \,\frac{d H(u)}{du}=\frac{r \delta u^{\delta-1}}{(1+u^\delta)^2}>0 \mbox{ and } \lim_{u\rightarrow\infty}H(u)=r>0,$$ we can apply Proposition 2.1 of Kang and Dieter (2011) to obtain the following conclusion: Assume that $u=H(u)$ has $n+1$ consecutive, distinct and non-degenerate solutions ${\bar{u}^i}, i=0,1,...,n$ with the following property {$$0=\bar{u}^0<\bar{u}^1<\cdot\cdot\cdot<\bar{u}^{n}.$$} then the even ${\bar{u}^{i}}, i\geq 2$ is locally asymptotical stable with $(\bar{u}^{i-1}, \bar{u}^{i+1})$ as its basin of attraction when ${\bar{u}^0=0}$ is locally asymptotical stable (i.e.,$\frac{d H(u)}{du}\big\vert_{u=0}<1$) and the odd ${\bar{u}^{i}}, i\geq 1$ is locally asymptotical stable with $(\bar{u}^{i-1}, \bar{u}^{i+1})$ as its basin of attraction when ${\bar{u}^0=0}$ is unstable (i.e.,$\frac{d H(u)}{du}\big\vert_{u=0}>1$).

According $\frac{d H(u)}{du}=\frac{r \delta u^{\delta-1}}{(1+u^\delta)^2}$, we can conclude that ${\bar{u}^0=0}$ is unstable if $\delta<1$ and ${\bar{u}^0=0}$ is locally asymptotical stable if $\delta>1$. Therefore, we have the following statement: i) If $\delta \in (0, 1)$ and Model \eqref{sA} has two non-negative equilibria $0, K$, then 0 is an unstable and
and $K$ is a stable positive equilibrium that attracts all points in $(0,\infty)$. ii) If $\delta>1$ and the only equilibrium of Model \eqref{sA} is 0, then $0$ is globally stable. iii) If  $\delta>1$ and Model \eqref{sA} has three non-negative distinct equilibria: $0,\, A$ and $K$ such that $0 < A < K$, then 0 equilibrium is locally asymptotical stable with the basin of attraction $[0, A)$; $A$ is unstable (repellor), while $K$ is locally asymptotical stable with the basin of attraction $(A ,\infty)$.

Now we only need to derive the conditions  when \eqref{sA} has only one, or  two or three equilibria. As we know, 0 is always a equilibrium of \eqref{sA}. If $u^*$ is a nontrivial positive equilibrium of \eqref{sA}, then we have $h(u^*)=1$, i.e., $u^*$ is a positive root of $F(u)=u^\delta-ru^{\delta-1}+1$. If $\delta<1$, then we have 
$$\lim_{u\rightarrow 0^+}F(u)=-\infty, \,\lim_{u\rightarrow \infty}F(u)=\infty \mbox{ and } \frac{dF}{du}=\delta u^{\delta-1}-r(\delta-1)u^{\delta-2}>0.$$ This implies that $F(u)=u$ has a unique positive solution $u^*$ and $$\frac{dH}{du}\big\vert_{u=u^*}=\frac{r \delta (u^*)^{\delta-1}}{(1+(u^*)^\delta)^2}=\frac{\delta}{1+(u^*)^\delta}<\delta<1.$$ Therefore, if $\delta \in (0, 1)$, Model \eqref{sA} has only two non-negative equilibria $0, K$ where 0 is an unstable and and $K$ is a stable positive equilibrium that attracts all points in $(0,\infty)$.

If $\delta>1$, then $\frac{dF}{du}=\delta u^{\delta-1}-r(\delta-1)u^{\delta-2}=u^{\delta-2}[\delta u-r(\delta-1)]=0$ when $u=u_c=\frac{r(\delta-1)}{\delta}$. Notice that 
$$F(0)=1>0, \, \frac{dF}{du}\big\vert_{u<u_c}<0 \mbox{ and }  \frac{dF}{du}\big\vert_{u>u_c}>0,$$ thus $F(u)=0$ (i.e., $h(u)=1$) have no positive root if 
$$F(u_c)>0\Rightarrow h(u_c)<1\Rightarrow \frac{r u_c^{\delta-1}}{1+u_c^\delta}=\frac{r  \left(\frac{r(\delta-1)}{\delta}\right)^{\delta-1}}{1+ \left(\frac{r(\delta-1)}{\delta}\right)^\delta}<1\Rightarrow r<r_{crit}=\delta (\delta-1)^{\frac{1}{\delta}-1}$$ and $F(u)=0$ have two non-negative distinct equilibria $0 < A<u_c < K$ if 
$$F(u_c)<0\Rightarrow h(u_c)>1\Rightarrow \frac{r u_c^{\delta-1}}{1+u_c^\delta}=\frac{r  \left(\frac{r(\delta-1)}{\delta}\right)^{\delta-1}}{1+ \left(\frac{r(\delta-1)}{\delta}\right)^\delta}>1\Rightarrow r>r_{crit}=\delta (\delta-1)^{\frac{1}{\delta}-1}.$$
Therefore the statements of Proposition \ref{th1:sA} hold.
\end{proof}

\subsection*{Proof of Theorem \ref{th2:sA-mss}}
\begin{proof}
First, we can easily verify that Model \eqref{sA-mss} is positively invariant and bounded in $\mathbb R_+$. If $u^*$ is a positive equilibrium of \eqref{sA-mss}, then it is a positive root of $F(u)=b u^{d}+u^\delta-r u^{\delta-1}+1$. Since
$$F'(u)=bdu^{d-1}+\delta u^{\delta-1}-r(\delta-1) u^{\delta-2}=u^{\delta-2}[bdu^{d-\delta+1}+\delta u-r(\delta-1)]=u^{\delta-2}G(u)$$where $G(u)=bdu^{d-\delta+1}+\delta u-r(\delta-1).$ If $\delta<1$, then we have
$$\lim_{u\rightarrow 0^+}F(u)=-\infty, \, F'(u)>0,\, \lim_{u\rightarrow 0^+}F'(u)=\infty \mbox{ and }\lim_{u\rightarrow \infty}F(u)=\infty,$$thus, $F(u)$ has a unique positive root $K>0$. 

Moreover, we have $H(0)=0$ and $H'(u)=\frac{r  u^{\delta-1}[\delta+b u^{d} (\delta-d)]}{(1+u^\delta+bu^{d})^2}$. This indicates the following two cases:
\begin{enumerate}
\item If $\delta\geq d$, then $H'>0$, thus, Model \eqref{sA-mss} is a monotone system where its dynamics can be classified into the following two cases:
\begin{enumerate}
\item If $\delta<1$, then Model \eqref{sA-mss} has two non-negative equilibria 0 and $K$ where 0 is unstable and $K$ is globally stable.
\item If $\delta>1+d$, then we have $u=0$ is locally stable and there exists some unique positive $u^c$ such that $G'(u)\vert_{u<u^c}<0 \mbox{ and } G'(u)\vert_{u>u^c}>0.$
Since $\lim_{u\rightarrow 0^+}G(u)=\infty$ and $\lim_{u\rightarrow \infty}G(u)=\infty$, therefore, $G(u)$ has exactly two positive roots if $G(u_c)<0$. Also, we can conclude that $G(u)>0$ if $G(u_c)>0$. This implies that if $G(u_c)<0$, then $F(u)$ has exactly two positive distinct critical points $0<u_c^1<u^c<u_c^2$. Therefore, we can conclude that 1) $F(u)$ has exact four distinct positive roots if $F(u_c^1)<0$ and 
$F(u_c^2)<0$; 2) $F(u)$ has exact two distinct positive roots if $F(u_c^1)F(u_c^2)<0$; and 3) $F(u)$ has no positive root if $F(u_c^1)>0$ and 
$F(u_c^2)>0$. In the case that $F(u)$ has no positive root, we can conclude that $u=0$ is globally stable. While if $G(u_c)>0$, then $F(u)<u$ for all $u>0$, therefore, we have $u=0$ is globally stable.
\item If $1+d>\delta>1$, then we have $u=0$ is locally stable and
$$G(0)=-r(\delta-1),\,G'(u)>0 \mbox{ and } \lim_{u\rightarrow \infty}G(u)=\infty,$$ thus, $F'(u)$ has a unique positive point $u_c$ such that $G(u_c)=0\Rightarrow F'(u_c)=0$ and 
$$F'(u)\vert_{u<u_c}>0 \mbox{ and } F'(u)\vert_{u>u_c}<0.$$Therefore, if $h(u_c)>1$, then Model \eqref{sA-mss} has as least two consecutive, distinct and non-degenerate equilibria $A,K$ such that $0<A<u_c<K$.
\end{enumerate}
\item If $\delta<d$, then there exists a unique critical point $u_c= (\frac{r(d-\delta)}{\delta})^{1/d}$ such that
$$H'(u)\vert_{u<u_c}>0, \,\, H'(u_c)=0\mbox{ and } H'(u)\vert_{u>u_c}<0.$$We consider the following two cases:
\begin{enumerate}
\item If $\delta<1$, then Model \eqref{sA-mss} is a unique hump model that is similar to the Ricker's map, thus, it has two non-negative equilibria 0 and $K$ where 0 is unstable.
\item If $\delta>1$, then $H'(0)=0$, thus $u=0$ is locally asymptotically stable. If  $r>r_{crit}$, then we have $H(u_c)>u_c$. This implies that Model \eqref{sA-mss} has at least two distinct positive roots $A,K$ such that $0<A<u_c<K$.
\end{enumerate}
\end{enumerate}

If both $\delta$ and $d$ are positive integers, then the number of interior equilibrium of \eqref{sA-mss} is the same as the number of positive roots of the following polynomial:
$$F(x)=bx^d+x^\delta-rx^{\delta-1}+1.$$
Descartes' rule of signs (Smith and Latham 1954; Meserve 1982) states that if the terms of a single-variable polynomial with real coefficients are ordered by descending variable exponent, then the number of positive roots of the polynomial is either equal to the number of sign differences between consecutive nonzero coefficients, or is less than it by a multiple of 2.  Therefore, the statement of Part 4 holds.

\end{proof}
\subsection*{Proof of Theorem \ref{th3:cs}}
\begin{proof}
Denote by $T: \mathbb R^2_+\rightarrow \mathbb R^2_+$ the map defined by the right side of Model \eqref{gx}-\eqref{gy}, then we have
\bae\label{T}
T&=&(T_1(x,y), T_2(x,y))=\left(\frac{r_1 x^{\delta_1}}{1+x^{\delta_1}+b_1y^{\delta_3}}, \frac{r_2 y^{\delta_2}}{1+y^{\delta_2}+b_2x^{\delta_4}}\right).
\eae  Thus, $T$ is $C^1$ and $T_i<r_i, i=1,2$. The Jacobian matrix of $T$ can be represented as follows:
\bae\label{Jg}
J&=&\left[\begin{array}{lcl}
\frac{r_1\delta_1x^{\delta_1-1}(1+b_1y^{\delta_3})}{\left(1+x^{\delta_1}+b_1y^{\delta_3}\right)^2}&-\frac{a_1b_1\delta_3x^{\delta_1}y^{\delta_3-1}}{\left(1+x^{\delta_1}+b_1y^{\delta_3}\right)^2}\\
-\frac{a_2b_2\delta_4y^{delta_2}x^{\delta_4-1}}{\left(1+y^{\delta_2}+b_2x^{\delta_4}\right)^2}&\frac{a_2\delta_2y^{\delta_2-1}(1+b_2x^{\delta_4})}{\left(1+y^{\delta_2}+b_2x^{\delta_4}\right)^2}
\end{array}
\right]
\eae which gives
\be\label{dJ}
det(J)=\frac{r_1r_2[\delta_1\delta_2(x^{\delta_1-1}y^{\delta_2-1}+b_2x^{\delta_1+\delta_4-1}y^{\delta_2-1}+b_1x^{\delta_1-1}y^{\delta_2+\delta_3-1})+b_1b_2x^{\delta_1+\delta_4-1}y^{\delta_2+\delta_3-1}(\delta_1\delta_2-\delta_3\delta_4)]}{\left(1+x^{\delta_1}+b_1y^{\delta_3}\right)^2\left(1+x^{\delta_1}+b_1y^{\delta_3}\right)^2}
\ee According to Lemma \ref{l1:bp}, we know that the map $T$ is positively invariant in $Int\,\mathbb R^2_+$, i.e., $T:Int\, \mathbb R^2_+\rightarrow Int \,\mathbb R^2_+$. Thus, for $(x,y)\in Int\,\mathbb R^2_+$, we have $det(J)\geq0$ (where equality occurs when $\delta_1\delta_2=\delta_3\delta_4$) and $J $ has a constant sign configuration as $J=\left[\begin{array}{lcl}
+&-\\
-&+
\end{array}
\right].$ This indicates that \eqref{Jg} is $K$-positive in $Int\,\mathbb R^2_+$ if and it is $K$-strongly positive if $\delta_1\delta_2>\delta_3\delta_4$.

Take $u,v\in Int\,\mathbb R^2_+$ such that $u<_Kv$, then by the
mean value theorem and \eqref{Jg}, we have
$$
\begin{array}{lcl}
T_1(v)-T_1(u)&=&\frac{\partial T_1}{\partial x}(v_1-u_1)+\frac{\partial T_1}{\partial y}(v_2-u_2)>0\\
T_2(v)-T_2(u)&=&\frac{\partial T_2}{\partial x}(v_1-u_1)+\frac{\partial T_2}{\partial y}(v_2-u_2)<0
\end{array}
$$where the partial derivatives are evaluated at a point on the line segment connecting $u$ and $v$. This implies that $T(u)<<_KT(v)$, i.e., the map $T$ is strongly competitive.

Choose any $(x_0,y_0)\in T(Int\, \mathbb R^2_+)$, we have $x_0<r_1$ and $y_0<r_2$. Let $(x,y)$ be a preimage of $(x_0, y_0)$, then it satisfies the following equation:
\bae\label{px}
x_0&=&\frac{r_1 x^{\delta_1}}{1+x^{\delta_1}+b_1y^{\delta_3}}\\
\label{py}
 y_0&=&\frac{r_2 y^{\delta_2}}{1+y^{\delta_2}+b_2x^{\delta_4}}
\eae This indicates that

\bae\label{rpx}
x_0&=&\frac{r_1 x^{\delta_1}}{1+x^{\delta_1}+b_1(\frac{y_0(1+b_2x^{\delta_4})}{r_2-y_0})^{\delta_3/\delta_2}}=f(x)\\
\label{rpy}
 y_0&=&\frac{r_2 y^{\delta_2}}{1+y^{\delta_2}+b_2(\frac{x_0(1+b_1y^{\delta_3})}{r_1-x_0})^{\delta_4/\delta_1}}=g(y).
\eae Since $\delta_1\delta_2\geq\delta_3\delta_4$, we have
\bae\label{dfg}\begin{array}{lcl}f'(x)&=&\frac{r_1(\delta_1\delta_2[x^{\delta_1}+b_2x^{\delta_1+\delta_4}+b_1 x^{\delta_1}(\frac{y_0(1+b_2x^{\delta_4})}{r_2-y_0})^{\delta_3/\delta_2}]+b_1b_2x^{\delta_1+\delta_4}(\frac{y_0(1+b_2x^{\delta_4})}{r_2-y_0})^{\delta_3/\delta_2}[\delta_1\delta_2-\delta_3\delta_4])}{\delta_2 x (1+b_2x^{\delta_4})(1+x^{\delta_1}+b_1(\frac{y_0(1+b_2x^{\delta_4})}{r_2-y_0})^{\delta_3/\delta_2})^2}>0\\
g'(y)&=&\frac{r_2(\delta_1\delta_2[y^{\delta_2}+b_1y^{\delta_2+\delta_3}+b_2 y^{\delta_2}(\frac{x_0(1+b_1y^{\delta_3})}{r_1-x_0})^{\delta_4/\delta_1}]+b_1b_2y^{\delta_2+\delta_3}(\frac{x_0(1+b_1y^{\delta_3})}{r_1-x_0})^{\delta_4/\delta_1}[\delta_1\delta_2-\delta_3\delta_4])}{\delta_1 y (1+b_1y^{\delta_3})(1+y^{\delta_2}+b_2(\frac{x_0(1+b_1y^{\delta_3})}{r_1-x_0})^{\delta_4/\delta_1})^2}>0.\end{array}\eae
This implies that both $f^{-1}(x_0)$ and $g^{-1}(y_0)$ have only one point, i.e., $(x_0,y_0)\in T(Int\, \mathbb R^2_+)$ has only one preimage. Thus, the map $T$ is injective in $Int\, \mathbb R^2_+$.

It is easy to check that $Int\, \mathbb R^2_+$ contains order intervals and is $\leq_K$-convex (i.e., it contains the line segment jointing any two of its points that are $K$-ordering). Therefore, we can apply Corollary 4.4 of Smith (1998) to obtain that $\{T^n(x_0,y_0)=(x_n, y_n)\}$ is eventually componentwise monotone for every $(x_0,y_0)\in Int\, \mathbb R^2_+$, i.e., there exists a positive integer $N$ such that either $x_n\leq x_{n+1}$ for all $n>N$ or $x_{n+1}\leq x_{n}$ for all $n>N$ and similarly for $y_n$. Notice that 
$$T(Int\, \mathbb R^2_+)\subset T( \mathbb R^2_+)\subset [0, r_1]\times[0, r_2],$$ thus, $\{T^n(x_0,y_0)=(x_n, y_n)\}$ converges to an equilibrium of Model \eqref{gx}-\eqref{gy} in $[0, r_1]\times[0, r_2]$. 
\end{proof}

\subsection*{Proof of Theorem \ref{th5:cs2}}
\begin{proof}
If $\delta_1\delta_2<\delta_3\delta_4$, then \eqref{dJ} vanishes at the following curves :
$$\begin{array}{lcl}
LC_{-1}^1:&& xy=0 \\
LC_{-1}^2:&&\delta_1\delta_2[1+b_2x^{\delta_4}+b_1y^{\delta_3}]+b_1b_2y^{\delta_3}x^{\delta_4}[\delta_2\delta_1-\delta_4\delta_3]=0.
\end{array}$$
Thus, we have
$$LC_{-1}^2:\, y=\gamma(x)=\left(\frac{\delta_1\delta_2[1+b_2x^{\delta_4}]}{b_1[b_2x^{\delta_4}(\delta_4\delta_3-\delta_2\delta_1)-\delta_1\delta_2]}\right)^{1/\delta_3},\mbox{ provided that } x>\left(\frac{\delta_1\delta_2}{b_2[\delta_3\delta_4-\delta_1\delta_2]}\right)^{1/\delta_4}$$ or
$$LC_{-1}^2:\, x=\gamma^{-1}(y)=\left(\frac{\delta_1\delta_2[1+b_1y^{\delta_3}]}{b_2[b_1y^{\delta_3}(\delta_4\delta_3-\delta_2\delta_1)-\delta_1\delta_2]}\right)^{1/\delta_4},\mbox{ provided that } y>\left(\frac{\delta_1\delta_2}{b_1[\delta_3\delta_4-\delta_1\delta_2]}\right)^{1/\delta_3}$$
 which gives that
$$\gamma'(x)=-\frac{b_2\delta_1\delta_2\delta_4^2x^{\delta_4-1}}{b_1[b_2x^{\delta_4}(\delta_4\delta_3-\delta_2\delta_1)-\delta_1\delta_2]^2}\left(\frac{\delta_1\delta_2[1+b_2x^{\delta_4}]}{b_1[b_2x^{\delta_4}(\delta_4\delta_3-\delta_2\delta_1)-\delta_1\delta_2]}\right)^{1/\delta_3-1}<0$$
 \mbox{ provided that } $x>\left(\frac{\delta_1\delta_2}{b_2[\delta_3\delta_4-\delta_1\delta_2]}\right)^{1/\delta_4}.$
In addition, by straight forward calculations, we have $$\gamma''(x)>0\mbox{ if }\delta_3\geq\delta_4\geq 1.$$ The curve $\gamma(x)$ separates $Int \,\mathbb R^2_+$ into two connected unbounded components. The one unbounded component containing $(0,\epsilon)\times (0,\epsilon)$ for very small $\epsilon$ is denoted by $\Omega_+$ and it has $det(J)>0$. While the other unbounded component is denoted by $\Omega_-$ and it has $det(J)<0$.

Let $T$ be the map describing Model \eqref{gx}-\eqref{gy}. The image of $\gamma(x)$ under the map $T$, i.e., $T(\gamma)$, called the critical curves in (Mira \emph{et al} 1996; Smith 1998; Kang 2011; Kang and Smith 2011), together with the portions of the coordinate axes which it cuts off, form the boundary of $T(Int\,\mathbb R^2_+)$. Because the curve $\gamma$ is linearly order by $<<_K$ (i.e., $\gamma'(x)<0$ indicates that any two points in $\gamma$ are related in the sense of strictly $K$-ordering) and $T$ preserves the strong ordering, $T(\gamma)$ is linearly ordered by $<<_K$ and connects $(0,r_2)$ to $(r_1,0)$. Let $B$ denote the bounded component of $Int\,\mathbb R^2_+\setminus T(\gamma)$. Then $T(\Omega_+)=T(\Omega_-)=B$, i.e., the map $T$ folds $Int\,\mathbb R^2_+$ over at $\gamma$ and maps both unbounded components onto $B$.

 If  $$r_1<\left(\frac{\delta_1\delta_2}{b_2[\delta_3\delta_4-\delta_1\delta_2]}\right)^{1/\delta_4}\mbox{ or }r_2<\left(\frac{\delta_1\delta_2}{b_1[\delta_3\delta_4-\delta_1\delta_2]}\right)^{1/\delta_3},$$ then $B$ is disjoint from $\gamma$ and is totally contained in $\Omega_+$. This implies that $det(T)_{B}>0$. Note that the bounded component of $Int\,\mathbb R^2_+\setminus T(\gamma)$ contains the point without pre-image (e.g., choosing a point with x-coordinate larger than $r_1$ and y-coordinate larger than $r_2$), then all the points in $Int\,\mathbb R^2_+\setminus T(\gamma)$ has zero pre-image. This because that the cardinality of $T^{-1}(w), w\in Int\,\mathbb R^2_+\setminus T(\gamma)$ is finite and constant according to Lemma 3.4 of Chow and Hale (1982). Similarly, we can show that the cardinality of $T^{-1}(w), w\in B$ is two where one is located in $B$ and the other one is located in $Int\,\mathbb R^2_+\setminus T(\gamma)$. This implies that
 $$T: B\rightarrow B \mbox{ is injective }.$$

 Therefore, we apply Corollary 4.4 of Smith (1998) to obtain that $\{T^n(x_0,y_0)=(x_n, y_n)\}$ is eventually componentwise monotone for every $(x_0,y_0)\in Int\, \mathbb R^2_+$. This indicates that every orbit of Model \eqref{gx}-\eqref{gy} with any initial condition in $\mathbb R^2_+$ converges to one of its equilibria. 
\end{proof}
\subsection*{Proof of Theorem \ref{th4:bs}}
\begin{proof}
For any initial value $(x_0,y_0)$ taken in $O_{ex}$, we have $x_0<A_1$, thus we have
$$x_1=T_1(x_0,y_0)=\frac{r_1 x_0^{\delta_1}}{1+x_0^{\delta_1}+b_1y_0^{\delta_3}}\leq T_1(x_0,0).$$
Since $\frac{\partial T_1}{\partial x}\vert_{x>0}>0$, thus we have
$$x_{t+1}=T_1^t(x_0,y_0)\leq T_1^t(x_0,0).$$
According to Proposition \ref{th1:sA}, we have $\lim_{t\rightarrow\infty}T_1^t(x_0,0)=0$. Since Model \eqref{gx}-\eqref{gy} is positively invariant in $\mathbb R^2_+$ by Lemma \ref{l1:bp}, thus, $\lim_{t\rightarrow\infty}x_t=0$ holds for Model \eqref{gx}-\eqref{gy} if Condition \textbf{H1} is satisfied and its initial condition is taken in $O_{ex}$.

Similarly, we can use the same argument to show that $\lim_{t\rightarrow\infty}y_t=0$ holds for Model \eqref{gx}-\eqref{gy} if Condition \textbf{H1} is satisfied and its initial condition is taken in $O_{ey}$. Therefore,
$$\lim_{t\rightarrow\infty}(x_t,y_t)=(0,0)=E_0 \mbox{ provided } (x_0,y_0)\in O_0.$$

Based on our arguments above, we can conclude that for any initial value $(x_0,y_0)$ taken in $O_{x}$, we have $\lim_{t\rightarrow\infty}y_t=0.$ In addition, we have
$$x_1=T_1(x_0,y_0)=\frac{r_1 x_0^{\delta_1}}{1+x_0^{\delta_1}+b_1y_0^{\delta_3}}\geq T_1(x_0, A_2)=\frac{r_1 x_0^{\delta_1}}{1+x_0^{\delta_1}+b_1(A_2)^{\delta_3}}.$$
Then by induction, we have $x_{t+1}=T_1^t(x_0,y_0)\geq T_1^t(x_0,A_2).$ Since Model \eqref{gx}-\eqref{gy} satisfies Condition \textbf{H2} and $x_0>A_1^{a_1}$, then according to Corollary \ref{c1:sA-a}, we have 
$$\lim_{t\rightarrow\infty}T_1^t(x_0,A_2)=K_1^{a_1}>A_1^{a_1}>A_1.$$This indicates that for Model \eqref{gx}-\eqref{gy}, we have
$$\liminf_{t\rightarrow\infty}x_t\geq K_1^{a_1}>A_1^{a_1}>A_1.$$Notice that $\lim_{t\rightarrow\infty}y_t=0$, thus, for any $\epsilon>0$, there exists a $N$ large enough, such that for any $t>N$, we have
$$x_{t+1}\geq \frac{r_1 x_t^{\delta_1}}{1+x_t^{\delta_1}+\epsilon}\geq A_1^{a_1}>A_1.$$
Let $\epsilon\rightarrow 0$, then the limiting system of Model \eqref{gx}-\eqref{gy} becomes $T_1(x,0)$. Since Condition \textbf{H1} is satisfied and the initial value is larger than $A_1$, thus according to Proposition \ref{th1:sA}, we have $\lim_{t\rightarrow\infty}x_t=K_1.$

Similarly, we can apply the same argument above to show the case when an initial value $(x_0,y_0)$ is taken in $O_{x}$. Therefore, the statement of Theorem \ref{th4:bs} holds.
\end{proof}
\subsection*{Proof of Theorem \ref{th5:bsE0}}
\begin{proof}
Assume that Model \eqref{gx}-\eqref{gy} satisfies Condition \textbf{H1}. Let $y=kx$, then from \eqref{T} any point $(x, kx)\in Int\,\mathbb R^2_+$ maps to the following point 
$$T(x,kx)=\left(T_1(x,kx),T_2(x,kx)\right)=\left(\frac{r_1 x^{\delta_1}}{1+x^{\delta_1}+b_1 (kx)^{\delta_3}},\frac{r_2 (kx)^{\delta_2}}{1+(kx)^{\delta_2}+b_2 x^{\delta_4}}\right).$$
If $\delta_3>\delta_1$, then we have
$$\frac{\partial T_1 (x,kx)}{\partial x}=\frac{r_1 x^{\delta_1-1}[\delta_1+b_1(kx)^{\delta_3}(\delta_1-\delta_3)]}{(1+x^{\delta_1}+b_1 (kx)^{\delta_3})^2}$$ which indicates that 
$$\frac{\partial T_1 (x,kx)}{\partial x}\vert_{x<x_c}>0\mbox{ and }\frac{\partial T_1 (x,kx)}{\partial x}\vert_{x>x_c}<0, x_c=\left(\frac{\delta_1}{b_1k^{\delta_3}(\delta_3-\delta_1)}\right)^{1/\delta_3}.$$
Let $h(x)=\frac{r_1 x^{\delta_1}}{x^{\delta_1}+b_1 (kx)^{\delta_3}}=\frac{r_1}{1+b_1k^{\delta_3}x^{\delta_3-\delta_1}}>T_1(x,kx)$.
Since $h(x)>T_1(x,kx)$ and $h'(x)<0$ provided $\delta_3>\delta_1$, thus we have
$$T_1(x, kx)<h(x)<A_1\Rightarrow x>\left(\frac{r_1-A_1}{b_1k^{\delta_3}}\right)^{\frac{1}{\delta_3-\delta_1}}.$$This gives 
$$T_1(x,kx)<A_1 \mbox{ whenever } x>\left(\frac{r_1-A_1}{b_1k^{\delta_3}}\right)^{\frac{1}{\delta_3-\delta_1}}.$$
Define $$O_{ex}^l=\bigcup_{k>0}\{(x,kx)\in Int\,\mathbb R^2_+: x>\left(\frac{r_1-A_1}{b_1k^{\delta_3}}\right)^{\frac{1}{\delta_3-\delta_1}}\}.$$
Then we can conclude that if Model \eqref{gx}-\eqref{gy} satisfies Condition \textbf{H1} and $\delta_3>\delta_1$, then for any initial value taken in $O_{ex}^l$, we have
$$\lim_{t\rightarrow\infty}x_t=0.$$
Similarly, we can conclude that if Model \eqref{gx}-\eqref{gy} satisfies Condition \textbf{H1} and $\delta_4>\delta_2$, then for any initial value taken in $O_{ey}^l$ where
$$O_{ey}^l=\bigcup_{k>0}\{(ky,y)\in Int\,\mathbb R^2_+: y>\left(\frac{r_2-A_2}{b_2k^{\delta_4}}\right)^{\frac{1}{\delta_4-\delta_2}}\},$$
 we have
$$\lim_{t\rightarrow\infty}y_t=0.$$
Let $O_0^l=O_{ex}^l\cap O_{ey}^l.$Then based on the arguments above, we can conclude that if Model \eqref{gx}-\eqref{gy} satisfies Condition \textbf{H1} and $\delta_3>\delta_1, \delta_4>\delta_2$, then for any initial value taken in $O_{0}^l$, we have
$$\lim_{t\rightarrow\infty}(x_t,y_t)=E_0.$$The argument above combined with the results of Theorem \ref{th4:bs}, we can conclude that the statement is hold.
\end{proof}
\subsection*{Proof of Theorem \ref{th6:ie}}
\begin{proof}Since Model \eqref{gx}-\eqref{gy} satisfies Condition \textbf{H1}, thus if $(x,y)$ is an interior equilibrium of \eqref{gx}-\eqref{gy}, we have the the following two equations:
$$\begin{array}{lcl}
1&=&\frac{r_1 x^{\delta_1-1}}{1+x^{\delta_1}+b_1y^{\delta_3}}\Rightarrow y=\left(\frac{r_1x^{\delta_1-1}-x^{\delta_1}-1}{b_1}\right)^{1/\delta_3}=F_1(x),\, 0<A_1<x<K_1\\
1&=& \frac{r_2 y^{\delta_2-1}}{1+y^{\delta_2}+b_2x^{\delta_4}}\Rightarrow x=\left(\frac{r_2y^{\delta_2-1}-y^{\delta_2}-1}{b_2}\right)^{1/\delta_4}=F_2(y),\, 0<A_2<y<K_2.\end{array}$$
Since
    $$\begin{array}{lcl}
    F_1'(x)& =&\frac{x^{\delta_1-2}\left(\frac{r_1x^{\delta_1-1}-x^{\delta_1}-1}{b_1}\right)^{1/\delta_3-1}}{b_1\delta_3}\left[r_1(\delta_1-1)-\delta_1 x\right]\\
     F_2'(y)& =&\frac{y^{\delta_2-2}\left(\frac{r_2y^{\delta_2-1}-y^{\delta_2}-1}{b_2}\right)^{1/\delta_4-1}}{b_2\delta_4}\left[r_2(\delta_2-1)-\delta_2 y\right],
     \end{array}$$
     thus we have, 
$$\frac{dF_1}{dx}\vert_{x<x_c}>0,\,\frac{dF_1}{dx}\vert_{x=x_c}=0 \mbox{ and }  \frac{dF_1}{dx}\vert_{x>x_c}<0,$$ and 
$$\frac{dF_2}{dy}\vert_{y<y_c}>0,\,\frac{dF_2}{dy}\vert_{y=y_c}=0 \mbox{ and }  \frac{dF_2}{dy}\vert_{y>y_c}<0$$where
$$x_c=\frac{r_1(\delta_1-1)}{\delta_1} \mbox{ and }y_c=\frac{r_2(\delta_2-1)}{\delta_2}.$$
Therefore, according to the geometry of $y=F_1(x)$ and $x=F_2(y)$, we have follows:
\begin{itemize}
\item If $F_1(x_c)<A_2$ or $F_2(y_c)<A_1$, then Model \eqref{gx}-\eqref{gy} has no interior equilibrium.
\item If $A_2<F_1(x_c)<K_2, \,K_1<F_2(y_c)$ or $K_2<F_1(x_c), \,A_1<F_2(y_c)<K_1$, then Model \eqref{gx}-\eqref{gy} has two interior equilibria.
\item If $F_1(x_c)>K_2$ and $F_2(y_c)>K_1$, then Model \eqref{gx}-\eqref{gy} has four interior equilibria.
\end{itemize}
\end{proof}

\subsection*{Proof of Theorem \ref{th7:sie}}
\begin{proof}Let $(x^*,x^*)$ be a symmetric interior equilibrium of Model \eqref{sgx}-\eqref{sgy}, then we have
$$r (x^*)^{\delta-1}=1+(x^*)^{\delta}+b(x^*)^d.$$
This indicates that the Jacobian matrix of Model \eqref{sgx}-\eqref{sgy} evaluated at the symmetric interior equilibrium $(x^*, x^*)$ can be represented as 
\bae\label{sJg}
J_{(x^*,x^*)}&=&\left[\begin{array}{lcl}
\frac{\delta(1+b(x^*)^{d})}{1+(x^*)^{\delta}+b(x^*)^{d}}&-\frac{bd(x^*)^{d}}{1+(x^*)^{\delta}+b(x^*)^{d}}\\
-\frac{bd(x^*)^{d}}{1+(x^*)^{\delta}+b(x^*)^{d}}&\frac{\delta(1+b(x^*)^{d})}{1+(x^*)^{\delta}+b(x^*)^{d}}
\end{array}
\right].\eae
The eigenvalues of \eqref{sJg} are
$$\lambda_1(x^*)=\frac{\delta+b (\delta-d) (x^*)^d}{1+(x^*)^{\delta}+b(x^*)^d}\mbox{ and } \lambda_2(x^*)=\frac{\delta+b (\delta+d) (x^*)^d}{1+(x^*)^{\delta}+b(x^*)^d}$$ where $\lambda_1(x^*)$ denotes the stability of $(x^*,x^*)$ on the invariant manifold $\Omega_{y=x}$ and the sign of $\lambda_2(x^*)-1$ denotes whether the eigenvector associated with $\lambda_2$ is pointing towards or away from $\Omega_{y=x}$.

Therefore, the statement of Theorem \ref{th7:sie} holds.
\end{proof}

\bibliographystyle{elsarticle-harv}

\begin{thebibliography}{00}
\bibitem{Ackleh2007}Ackleh A.S., Allen L.J.S and Carter J., 2007. Establishing a beachhead: A stochastic population model with an Allee effect applied to species invasion. \emph{Theoretical Population Biology}, \textbf{71}, 290-300.
\bibitem{Allee1949}Allee W.C., Emerson A.E., Park O., Park T. and Schmidt K.P., 1949. \emph{Principles of animal ecology}, Saunders (W.B.) Co Ltd, 1st Edition.
\bibitem{Amarasekare1998a}Amarasekare P., 1998a. Allee effects in metapopulation dynamics. \emph{The American Naturalist}, \textbf{152} 298--302.
\bibitem{Amarasekare1998b}Amarasekare P., 1998b. Interactions between local dynamics and dispersal: insights from single species models. \emph{Theoretical Population Biology}, \textbf{53},  44--59.
\bibitem{Begon1996}Begon M., Harper J.L. and Townsend C.R., 1996. \emph{Ecology: Individuals, Populations and Communities}. Blackwell Science Ltd., Oxford.
\bibitem{Bengtsson1989}Bengtsson J., 1989. Interspecific competition increases local extinction rate in a metapopulation system. \emph{Nature}, \textbf{340}, 713-715.
\bibitem{Cai2014}Cai Yongli, Banerjee M., Kang Y. and  Wang W., 2014. Spatiotemporal complexity in a predator-prey model with weak Allee effects.  \emph{Mathematical Biosciences and Engineering}, \textbf{11}(6), 1247-1274. 
\bibitem{Calow1998}Calow P., Falk D.A., Grace J. and Moore, P.D., 1998. \emph{The encyclopedia of ecology and environmental management}. Blackwell Science, Oxford.
\bibitem{Chesson1989}\modifyc{Chesson P. L. and Ellner S., 1989. Invasibility and stochastic boundedness in monotonic competition models, \emph{Journal of Mathematical Biology}, \textbf{27}, 117-138.}
\bibitem{Chow1982}Chow S.N. and Hale J.K., 1982. \emph{Methods of Bifurcation Theory}, Springer-Verlag, New York.
\bibitem{Clark2003}Clark D., Kulenovic  M.R.S. and Selgrade J.F., 2003. Global asymptotic behavior of a two-dimensional difference equation modelling competition, \emph{Nonlinear Analysis}, \textbf{52}, 1765-1776.

\bibitem{Cushing2004}Cushing J., Levarge S., Chitnis N. and Henson S.M., 2004. Some discrete competition models and competitive exclusion principle, \emph{Journal of Difference Equations and Applications}, \textbf{10}, 1139-1151.
\bibitem{Courchamp2009}Courchamp F., Berec L.and Gascoigne J., 2009. \emph{Allee effects in ecology and conservation}. Oxford University Press.
\bibitem{Cushing2007}Cushing J., Henson S.M. and Blackburn C.C., 2007. Multiple mixed-type attractors in a competition model. \emph{Journal of Biological Dynamics}, \textbf{1}, 347-362.
\bibitem{Cushing2012}Cushing J. and Hudson J., 2012. Evolutionary dynamics and strong Allee effects,  \emph{Journal of Biological Dynamics}, \textbf{6}(2), 941-958.
%
\bibitem{Cushing2015}Cushing J., 2015.  The evolutionary dynamics of a population model with a strong Allee effect, \emph{Mathematical Biosciences and Engineering}, \textbf{12}(4), 643-660.

\bibitem{Dancer1991}Dancer  E.N. and Hess P., 1991. Stability of fixed points for order-preserving discrete-time dynamical systems, \emph{J. Reine. Angew. Math.}, \textbf{419}, 125-139.
\bibitem{Dennis1989}Dennis B., 1989. Allee effects: population growth, critical density, and the chance of extinction. Nat. \emph{Res. Model.}, \textbf{3}, 481-538.
\bibitem{Dennis2002}Dennis B., 2002. Allee effects in stochastic populations. \emph{Oikos}, \textbf{96}, 389-401. 

\bibitem{Drake2004}Drake J.M., 2004. Allee effects and the risk of biological invasion. \emph{Risk Analysis}, \textbf{24} 795--802.

\bibitem{Egami2010}Egami C, 2010. Permanence of delay competitive systems with weak Allee effects.
\emph{ Nonlinear Analysis: Real World Applications}, \textbf{11}, 3936-3945.
\bibitem{EKL2005}Edelstein-Keshet L., 2005. \emph{Mathematical models in biology}, SIAM, Philadelphia.
\bibitem{Elaydi2010}Elaydi S.N. and Sacker R.J., 2010. Population models with Allee effect: A new model. \emph{Journal of Biological Dynamics}, \textbf{4}, 397-408.

\bibitem{Etienne2002}Etiemme R., Werthei B., Hemerik L. Schneider and Powell J., 2002. The interaction between dispersal, the Allee effect and scramble competition affects population dynamics. \emph{Ecological Modeling}, \textbf{148}, 153-168. 
\bibitem{Fagan2002}Fagan W.F., Lewis M.A., Neubert M.G. and van den Driessche P., 2002. Invasion theory and biological control. \emph{Ecology Letters}, \textbf{5}, 148-157. 
\bibitem{Greene2001}Greene C. and Stamps J.A., 2001. Habitat selection at low population densities. \emph{Ecology}, \textbf{82}, 2091-2100.

\bibitem{Gyllenberg1999}Gyllenberg M., Hemminki J.and Tammaru T., 1999. Allee effects can both conserve and create spatial heterogeneity in population densities. \emph{Theoretical Population Biology}, \textbf{56}, 231-242.
\bibitem{Hopf1993}\modifyc{Hopf F. A., Valone T. J. and Brown J. H., 1993. Competition theory and the structure of ecological communities, \emph{Evolutionary Ecology}, \textbf{7}, 142-154.}
\bibitem{Franke1991a}Franke J. and Yakubu A-A., 1991a. Global attractors in competitive {systems}. \emph{Nonlinear Analysis, Theory, Methods \& Application}, \textbf{16}, 111-129.

\bibitem{Franke1991b}Franke J. and Yakubu A-A., 1991b. Mutual exclusion versus coexistence for discrete competitive systems. \emph{Journal of Mathematical Biology}, \textbf{30}, 161-168.
\bibitem{Hassell1976}Hassell M.P., Lawton J.H. and  Beddington J.R., 1976. The components of arthropod predation: I. The prey death-rate. \emph{J. Anim. Ecol.}, \textbf{45}, 135-164.
\bibitem{Harry2012}Harry A.J., Kent C.M. and Kocic V.L., 2012. Global behavior of solutions of a periodically forced Sigmoid Beverton-Holt model, \emph{Journal of Biological Dynamics}, \textbf{6}, 212-234.
\bibitem{Hess1991}Hess P. and Lazer A.C., 1991. On an abstract competition model and applications, \emph{Nonlinear Anal: TMA}, \textbf{16}, 917-940.
\bibitem{Jang2006}Jang S.R.J., 2006. Allee effects in a discrete-time host-parasitoid model. \emph{Journal of Difference Equations and Applications}, \textbf{12}, 165-181.
\bibitem{Gascoigne2004} Gascoigne J. and Lipcius R.N., 2004. Allee effects in marine systems, \emph{Mar. Ecol. Prog. Ser.}, \textbf{269}, 49-59.
\bibitem{Jiang1987}Jiang H. and Rogers T.D., 1987. The discrete dynamics of symmetric competition in the plane. \emph{J. Math. Biol.}, \textbf{25}, 573-596.
\bibitem{Kang2010}Kang Y. and Chesson P., 2010. Relative nonlinearity and permanence. \emph{Theoretical Population Biology}, \textbf{78}, 26-35.
\bibitem{Kang2011}Kang Y. and Lanchier N., 2011. Expansion or extinction: deterministic and stochastic two-patch models with Allee effects. \emph{Journal of Mathematical Biology}, \textbf{62}, 925-973.
\bibitem{Kang2011a}Kang Y. and Armbruster D., 2011. Dispersal effects on a two-patch discrete model for plant-herbivore interactions. \emph{Journal of Theoretical Biology}, \textbf{268}, 84-97.
\bibitem{Kang2011c}Kang Y. and Yakubu Abdul-Aziz, 2011. Weak Allee effects and species coexistence, \emph{Nonlinear Analysis: Real World Applications}, \textbf{12}, 3329-3345.
\bibitem{Kang2012d}Kang Y. and Castillo-Chavez C., 2012.  Multiscale analysis of compartment models with dispersal, \emph{Journal of Biological Dynamics}, \textbf{6(2)}, 50-79.

\bibitem{Kang20123}Kang Y., 2013. Permanence of a general discrete-time two-species interaction model with nonlinear per-capita growth rates. \emph{Discrete and Continuous Dynamical Systems} - Series B, \textbf{18}(8), 2123-2142.
\bibitem{Kang2013}Kang Y., 2013. Scramble competitions can rescue endangered species subject to strong Allee effects, \emph{Mathematics Biosciences}, \textbf{241}(1), 75-87.

\bibitem{Kang2014f}Kang Y., Bhowmick R. A.,  Sasmal K. S. and Chattopadhyay J., 2014a. Host-Parasitoid systems with predation-driven Allee effects in host population.  \emph{Journal of Biological Dynamics}. [Epub ahead of print]. DOI:10.1080\/17513758.2014.972473

\bibitem{Kang2014e}Kang Y. and Udiani O., 2014. Dynamics of a single species evolutionary model with Allee effects. \emph{Journal of Mathematical Analysis and Applications}, \textbf{418} (1),  492-515.
\bibitem{Kang2014d}Kang Y., Sasmal K. S., Bhowmick R. A. and Chattopadhyay J., 2014b. Dynamics of a predator-prey system with prey subject to Allee effects and disease, \emph{Journal of Mathematical Biosciences and Engineering}, \textbf{11}(4), 877-918.

\bibitem{Kang2014a}Kang Y. and Carlos Castillo-Chavez C., 2014a. Dynamics of SI models with both horizontal and vertical transmissions as well as Allee effects. \emph{Journal of Mathematical Biosciences}, \textbf{248}, 97-116.
\bibitem{Kang2014b}Kang Y. and Carlos Castillo-Chavez C., 2014b. A simple epidemiological model for populations in the wild with Allee effects and disease-modified fitness. \emph{Discrete and Continuous Dynamical Systems} - Series B, \textbf{19} (1), 89-130.

\bibitem{Kang2014c}Kang Y. and Carlos Castillo-Chavez C., 2014c. A simple two-patch epidemiological model with Allee effects and disease-modified fitness. \emph{Mathematics of Continuous and Discrete Dynamical Systems}, edit by Abba B. Gumel in Contemporary Mathematics, 618.
%
\bibitem{Keitt2001}Keitt T.H., Lewis M.A. and Holt R.D., 2001. Allee effects, invasion pinning, and species' borders. \emph{The American Naturalist}, \textbf{157}, 203-216.

\bibitem{Liebhold2003}Liebhold A. and Bascompte J., 2003. The Allee effect, stochastic dynamics and the eradication of alien species. \emph{Ecology Letters}, \textbf{6}, 133-140.
\bibitem{Lidicker2010}Lidicker W.Z., 2010. The Allee effect: its history and future importance, \emph{The Open Ecology Journal}, \textbf{3}, 71-82.
\bibitem{Livadiotis2012}Livadiotis G. and Elaydi S., 2012. General Allee effect in two-species population biology. \emph{Journal of Biological Dynamics}, \textbf{6}(2), 959-973.
\bibitem{Livadiotis2015}Livadiotis G., Assas L., Dennis B., Elaydid S. and Kwessi E., 2015 discrete-time host-parasitoid model with an Allee effect.  \emph{Journal of Biological Dynamics}, \textbf{9}(1). DOI:10.1080\/17513758.2014.982219

\bibitem{Luis2009}Luis R., Elaydi S. and Oliveira H., 2009. Nonautonomous periodic systems with Allee effects. \emph{Journal of Difference Equations and Applications}, \textbf{1}, 1-16.

\bibitem{McCathy1997}McCarthy  M.A., 1997. The Allee effect, finding mates and theoretical models. \emph{Ecol. Model.}, \textbf{103}, 99-102.
\bibitem{Meserve1982}Meserve B. E., 1982. \emph{Fundamental Concepts of Algebra}, New York, Dover Publications.
\bibitem{Mira1996}Mira C., Gardini, Barugola A. and Cathala J.-C., 1996. \emph{Chaotic Dynamics in Two-Dimensional Non-invertible Maps}, Nonlinear Sciences, Series A., Vol. 20, World Scientific Publishing Co., Pte. Ltd, Singapore.
\bibitem{Myers1995}Myers R.A. , Barrowman N.J. , Hutchings J.A. and Rosenberg A.A., 1995. Population dynamics of exploited fish stocks at low population levels, \emph{Science}, \textbf{269}, 1106-1108.
\bibitem{Myers1998}Myers J.A. and Mertz G., 1998. Reducing uncertainty in the biological basis of fisheries management by meta-analysis of data from many population: a synthesis, \emph{Fish. Res.}, \textbf{37}, 51-60.
\bibitem{Myers2001}Myers R.A., 2001. Stock and recruitment: generalizations about maximum reproductive rate, density dependence, and variability using meta-analytic approaches, \emph{ICES J. Mar. Sci.}, \textbf{58}, 937-951.
\bibitem{Peng2015}Peng Feng and  Kang Yun, 2015. Dynamics of a modified Leslie-Gower model with double Allee effects. \emph{Nonlinear Dynamics}, \textbf{80}, (1-2), 1051-1062.
\bibitem{Petrovskii2005}Petrovskii S., Morozov A.and Li B.-L., 2005. Regimes of biological invasion in a predator-prey system with the Allee effect. \emph{Bulletin of Mathematical Biology}, \textbf{67}, 637-661.

%
\bibitem{Scheuring1999}Scheuring  I., 1999. Allee effect increases dynamical stability in populations. \emph{J. Theor. Biol.}, \textbf{199}, 407-414.

\bibitem{Schreiber2003}Schreiber S., 2003. Allee effects, extinctions, and chaotic transients in simple population models. \emph{Theoretical Population Biology}, \textbf{64}, 201-209.
\bibitem{Selgrade1987}Selgrade J.F. and Ziehe M., 1987. Convergence to equilibrium in a genetic model with differential viability between the sexes, \emph{J. Math. Biol.}, \textbf{25}, 477-490.
\bibitem{Selgrade1992}Selgrade J.F. and Namkoong G., 1992. Dynamical behavior for population genetics models of differential and difference equations with nonmonotone fitnesses. \emph{Journal of Mathematical Biology}, \textbf{30}, 815-826.

\bibitem{Shigesada1997}Shigesada N. and Kawasaki K., 1997. Introduction. In: N. Shigesada and K. Kawasaki, Editors, \emph{Biological Invasions: Theory and Practice}, Oxford University Press, New York, USA, 1-5.

\bibitem{Smith-Latham1954}Smith D. E. and Latham M.L., 1954. \emph{The Geometry of Rene Descartes with a facsimile of the first edition}, translated by Smith D.E. and Latham M.L., New York, Dover Publications.
\bibitem{Smith1995}Smith H.L., 1995. \emph{Monotone Dynamical Systems: An Introduction to the Theory of Competitive and Cooperative Systems}, American Mathematical Society, Providence, RI.
\bibitem{Smith1998}Smith H.L., 1998. Planar competitive and cooperative difference equations, \emph{Journal of Difference Equations and Applications}, \textbf{3}, 335-357.

\bibitem{Stephens1999}Stephens P.A., Sutherland W.J. and R.P. Freckleton R.P., 1999. What is the Allee effect? \emph{Oikos}, \textbf{87}, 185-190.
\bibitem{Stoner2000}Stoner A.W. and Ray-Culp M., 2000 Evidence for Allee effects in an over-harvested marine gastropod: density-dependent mating and egg production, \emph{Mar. Ecol. Prog. Ser.}, \textbf{202}, 297-302.
\bibitem{Taylor2005}Taylor C.M. and Hastings A., 2005. Allee effects in biological invasions. \emph{Ecology Letters}, \textbf{8}, 895-908.

\bibitem{Terescak1996}I. Teresc\'ak, 1996. Dynamics of $C^1$ smooth strongly monotone discrete-time dynamical systems, preprint. 
\bibitem{Thieme2003}Thieme H., 2003. \emph{Mathematics in Population Biology}, Princeton University Press, Princeton.
\bibitem{Thieme2009}Thieme H., Dhirasakdanon T, Han Z. and Trevino R., 2009. Species decline and extinction: synergy of infectious disease and Allee effect? \emph{Journal of Biological Dynamics}, \textbf{3}, 305-323.
\bibitem{Thomson1993}Thomson G.G., 1993. A proposal for a threshold stock size and maximum fishing mortality rate, in Risk Evaluation and
Biological Reference Points for Fisheries Management, \emph{Canad. Spec. Publ. Fish. Aquat. Sci.}, Vol. \textbf{120}, S.J. Smith,
J.J. Hunt, and D. Rivard, eds., 303-320.
\bibitem{Valone1995}Valone T.J. and Brown J.H., 1995. Effects of competition, colonization and extinction on rodent species diversity. \emph{Science}, \textbf{267}, 880-883.
\bibitem{Vandermeer2002}Vandermeer J., Evans M.A., Foster P., H\"o\"ok T., Reiskind M. and Wund M., 2002. Increased competition may promote species coexistence. \emph{PNAS}, \textbf{99}, 8731-8736.
\bibitem{Wang2001}Wang Y. and Jiang Jifa, 2001. The general properties of discrete-time competitive dynamical systems, \emph{Journal of Differential Equations}, \textbf{176}, 470-493.
\bibitem{Wang2002}Wang M.H, Kot M. and M. G. Neubert M.G., 2002. Integrodifference equations, Allee effects, and invasions. \emph{Journal of Mathematical Biology}, \textbf{44}, 150-168.

\bibitem{Zhou2004}Zhou S.R., Liu C.Z. and Wang G., 2004. The competitive dynamics of metapopulation subject to the Allee-like effect. \emph{Theoretical Population Biology}, \textbf{65}, 29-37. 
\end{thebibliography}



\end{document}